\documentclass[11pt,twoside]{article}
\usepackage{amsmath,amsfonts,amssymb,cases,mathdots,color,colordvi,url,tikz,adjustbox}
\usepackage{hyperref}
\usepackage{bm}
\usepackage{pgfpages}
\usepackage{graphics}
\usepackage{epsfig}
\usepackage{algorithm}
\usepackage{algorithmic}
\usepackage{booktabs,multirow}
\usepackage{caption, subcaption}
\usepackage{authblk}
\usepackage{ulem}
\usepackage{comment}

\usepackage[UTF8]{ctex}
\usepackage{svg}

\usepackage{fontspec}
\usepackage[authoryear,round]{natbib}

\usepackage{graphicx}
\usepackage{threeparttable}
\usepackage{adjustbox}
\usepackage{array}

\newcommand{\bes}{\begin{equation*}}
\newcommand{\ees}{\end{equation*}}
\newcommand{\be}{\begin{equation}}
\newcommand{\ee}{\end{equation}}
\newcommand{\ba}{\begin{eqnarray}}
\newcommand{\ea}{\end{eqnarray}}
\newcommand{\bal}{\begin{aligned}}
\newcommand{\eal}{\end{aligned}}
\newcommand{\bas}{\begin{eqnarray*}}
\newcommand{\eas}{\end{eqnarray*}}

\graphicspath{ {./Figures/} }

\def\L{{\mathcal L}}

\def\W{{\mathcal W}}

\def\R{\mathbb{R}}

\def\bfr{{\bf r}}
\def\bff{{\bf f}}

\def\bfa{{\bf a}}

\def\bfy{{\bf y}}
\def\bfz{{\bf z}}
\def\bfw{{\bf w}}

\def\bfmu{{\bf {\boldsymbol \mu}}}

\def\bfone{{\bf 1}}

\title{Decision-Focused Covariance Hybrid for Portfolio Optimization}
    \author[1]{Yingxiao Wang }{}
    \author[2]{Ranran Wang\thanks{Corresponding author}}{}
    \affil[1]{Academy for Advanced Interdisciplinary Studies, Peking University, Beijing 100971, China.}
    \affil[2]{National Engineerig Laboratory for Big Data Analysis and Applications, Peking University, Beijing 100971, China.}
\date{}

\begin{document}
	
	\maketitle

\thispagestyle{empty}

\begin{abstract} 
Covariance matrix estimation is central to mean-variance portfolio optimization, as error in the covariance matrix can substantially affect the portfolio decisions and out-of-sample performance. Structured covariance estimators and decision-focused data-driven models provide different approaches to constructing covariance matrices, but the direct hybrid of these two matrices has not been sufficiently studied. We propose a decision-focused covariance hybrid framework that combines a structured covariance estimate with a data-driven covariance matrix. The data-driven covariance is learned through a portfolio-level objective. We further extend the framework to an adaptive hybrid scheme. Empirical results show that the two component matrices provide distinct risk representations. The hybrid covariance matrix improves the out-of-sample Sharpe ratio. The adaptive hybrid scheme also achieves competitive performance across investment styles. Spectral analysis further shows that the hybrid matrix exhibits improved numerical conditioning. These findings highlight the potential of combining structured covariance information with decision-focused risk representations for portfolio construction.

\vspace{3mm}
\end{abstract}
{}
\section{Introduction}\label{introduction}

In mean-variance (MV) portfolio optimization \citep{markowitz1952}, the covariance matrix is a key input that determines how portfolio risk is represented and penalized. Errors in covariance estimates can affect portfolio weights, risk exposures, turnover, and out-of-sample performance \citep{broadie1993computing,kan2007optimal}. These effects can be particularly pronounced in high-dimensional situations, because covariance estimation becomes more difficult and errors may be amplified by the optimizer \citep{demiguel2009a,ardia2017the}. Therefore, the covariance input is important not only as a statistical estimate but also for the performance of the portfolios constructed using it. Understanding these decision consequences is important for constructing covariance inputs that support portfolio optimization.

Structured and data-driven models offer different ways to construct covariance inputs for portfolio optimization. Structured covariance models impose assumptions on the risk representation to reduce estimation complexity or improve estimation stability. For example, factor models represent asset risk through a lower-dimensional set of common factors and provide interpretable covariance inputs in high-dimensional settings \citep{fama1993common,carhart1997persistence,Fama2015five,de2021factor}. Data-driven methods, meanwhile, use historical returns and other available features to construct inputs. They can estimate the covariance matrix directly or learn components of a structured model, potentially allowing more flexible representations of covariance \citep{kim2018forecasting,dorta2018structured,kelly2019characteristics,monteiro2020stochastic,Gu2021autoencoder,gopal2024neuralfactors,zhu2024dynamic}. These approaches have provided useful tools for modeling the asset covariance and constructing portfolios.

However, the structured covariance models rely on specific structural assumptions, which may limit their ability to capture risk patterns that deviate from those assumptions. Data-driven models may provide a different representation, but their portfolio implications depend on how the models are constructed and trained. Existing approaches include direct data-driven covariance estimation, data-driven modeling of components of structured models, and covariance shrinkage or cleaning. These developments naturally raise a question: how does directly combining a structured covariance estimate with a data-driven covariance matrix affect the risk penalties assigned to portfolios and, consequently, portfolio allocations and out-of-sample performance?

To address this question, the hybrid covariance should be considered in the context of the portfolio decision it supports. The decision-focused learning (DFL) framework provides a way to incorporate the portfolio decision objectives into model training, rather than focusing solely on improving prediction accuracy \citep{anis2025end,kim2025estimating,yang2026neural,bongiorno2026end}. Several related studies on covariance estimation focus on the global minimum-variance (GMV) portfolio construction. However, investors typically seek to trade off return and risk rather than minimize risk alone. The mean-variance (MV) model provides a complementary setting in which the covariance input interacts with expected return predictions and the investor's risk aversion.

Against this background, we investigate the direct matrix-level hybrid of a structured covariance matrix and a decision-focused data-driven covariance matrix for MV portfolio optimization. To the best of our knowledge, we are not aware of prior work that directly combines these two types of covariance matrices at the matrix level in a DFL portfolio optimization framework. Our focus is not on matrix combination alone, but on a specific construction in which the data-driven covariance matrix is learned using a portfolio-level objective. We examine how the component matrices and their hybrid assign risk penalties to portfolios and how these differences affect portfolio allocations and out-of-sample performance.

The main contributions of this paper are threefold.

(i) We propose a decision-focused covariance hybrid framework for MV portfolio optimization. The framework directly integrates a structured covariance estimate and a data-driven covariance matrix at the matrix level, with the latter learned through a portfolio-level objective. Then, the hybrid matrix is input in the portfolio optimization layer. In addition, the framework is designed to accommodate alternative covariance models and optimization layers.

(ii) We empirically analyze the risk representations and portfolio implications of the two component covariance matrices and their hybrid. Using an FF3 covariance model and an iTransformer-based covariance model as a concrete implementation, we examine their covariance structures, the risk penalties they assign to portfolios, and their effects on portfolio allocations and out-of-sample performance. Spectral analysis is also conducted to investigate their numerical properties. In the main out-of-sample comparisons, the hybrid strategies achieve a higher Sharpe ratio than the traditional and data-driven benchmarks considered. During the two selected market downturns, they also exhibit smaller losses and lower maximum drawdowns than the external data-driven benchmarks. These analyses provide empirical insights into how the proposed hybrid covariance matrix influences portfolio decisions and performance.

(iii) We extend the framework with an adaptive hybrid mechanism that adjusts the hybrid weight according to the decision features. A two-stage training procedure is adopted: the return and data-driven covariance prediction modules are trained first; their parameters are then fixed while the adaptive hybrid mechanism is trained. Numerical experiments show that the adaptive approach delivers competitive out-of-sample performance in our settings. It provides a flexible alternative to fixed hybrid weights.

The remainder of this paper is organized as follows. Section~\ref{sec:background} reviews the background on portfolio optimization and related work on structured and data-driven covariance estimation models. Section~\ref{sec:methodology} presents the proposed framework and its implementation. Section~\ref{sec:empirical-study} describes the experimental design and discusses the results. Section~\ref{sec:conclusion} concludes this paper.

\section{Background}\label{sec:background}

\subsection{Mean-variance portfolio optimization}\label{sec:mvo}

Consider a market with $N$ risky assets. Let $\bfmu_t \in \R^N$ and $\Sigma_t \in \R^{N \times N}$ denote the estimations of the expected return vector and the covariance matrix of returns at time $t$, respectively. Denote the portfolio weight vector as $\bfw_t \in \R^N$. The classical MV framework determines the portfolio weights by balancing the expected return against the portfolio risk. In this paper, we consider the following long-only MV optimization problem:
\begin{equation}\label{eq:mv}
\begin{aligned}
    \min_{\bfw_t} \quad & -\bfmu_t^\top \bfw_t 
    + \lambda \bfw_t^\top \Sigma_t \bfw_t \\
    \text{s.t.} \quad & \bfone_N^\top \bfw_t = 1, \quad \bfw_t \geq \bf0 ,
\end{aligned}
\end{equation}
where $\lambda > 0$ is the risk aversion parameter. A larger $\lambda$ places greater emphasis on the risk penalty. In this paper, we assume the covariance matrices are symmetric and positive semidefinite. 

In Problem \eqref{eq:mv}, the return and covariance inputs jointly determine the objective assigned to each feasible portfolio. Estimation errors of the two inputs can lead to substantial changes in portfolio weights and out-of-sample performance \citep{broadie1993computing,kan2007optimal}. These effects are especially relevant in high-dimensional situations \citep{demiguel2009a,ardia2017the}. Consequently, evaluating a covariance input for portfolio construction requires attention not only to its statistical properties but also to its effect on the resulting investment decisions.

\subsection{Structured covariance model}\label{sec:structured_covariance_model}

Structured covariance models reduce the difficulty of covariance estimation by imposing restrictions on the covariance matrix or on its components. Common approaches include factor structures, shrinkage estimators, and sparse structures. Factor models describe common variation in asset returns through a relatively small number of factors, while shrinkage estimators combine a sample covariance matrix with a structured target to improve estimation properties. Sparse covariance approaches, in turn, impose restrictions on the dependence structure among assets. These methods provide different ways to balance model flexibility and estimation stability in high-dimensional situations \citep{fama1993common,ledoit2004honey,fan2008high,Fan2012Vast,ledoit2012nonlinear}.

Factor models are widely used in covariance estimation \citep{fama1993common,fan2008high,Fan2012Vast,de2021factor}. Consider a linear factor model with $N$ assets and $P$ factors, where $P \ll N$. Let $\bfr_t \in \R^N$ and $\bff_t \in \R^P$ denote the asset return vector and factor return vector, respectively. The model can be written as
\begin{align*}
    \mathbf{r}_t=\mathbf{a}+B\mathbf{f}_t+\boldsymbol{\varepsilon},
\end{align*}
where $\bfa \in \R^N$ is the intercept vector, $B \in \R^{N \times P}$ is the factor loading matrix, and $\boldsymbol{\varepsilon}$ is the vector of idiosyncratic returns. Assuming that the factor returns and idiosyncratic returns are uncorrelated, the covariance matrix implied by the model is
\begin{align*}
    \Sigma^{\mathrm{factor}} = B \Sigma^f B^\top + \Psi,
\end{align*}
where $\Sigma^f \in \R^{P \times P}$ is the factor covariance matrix and $\Psi \in \R^{N \times N}$ is the covariance matrix of idiosyncratic returns. The latter is often assumed to be diagonal or sparse.

A factor model represents the return covariance matrix as a function of the factor loading matrix, the factor return covariance matrix, and the idiosyncratic return covariance matrix. It can improve covariance estimation in the high-dimensional situations \citep{fan2008high,Fan2012Vast}. Observable economic factors can also help interpret common sources of portfolio risk. For example, the Fama-French three factors model \citep{fama1993common} uses market, size, and value factors to describe systematic variation in returns. In addition to observable factors, latent factors can be constructed through dimension-reduction methods \citep{laloux2000random,zou2006sparse,lopez2016building}. \citet{conlon2025asset} investigates the factor-based covariance matrices using methods including principal component analysis, sparse principal component analysis, and autoencoders. They find that the above methods also can improve covariance forecasting and minimum-variance portfolio performance in large-scale asset universes.

The effectiveness of a structured covariance model depends on whether its assumptions adequately represent the relevant risk patterns. For example, a factor model depends on the selected factors, the specification of factor loadings, and assumptions about idiosyncratic risks. A given structure may therefore be insufficient to capture correlations that vary with market conditions. Matrix combination provides a way to incorporate different sources of risk information. In conventional linear shrinkage estimation, the sample covariance matrix is combined with a structured target to improve estimation properties \citep{ledoit2004honey}. Such targets are often specified using relatively simple structures, such as the identity matrix. Recent advances in data-driven methods offer new ways to learn nonlinear correlations. Our approach is related to shrinkage at the matrix-combination level, but differs in the construction of the data-driven component and the objective used to learn it. Specifically, we consider a data-driven covariance matrix generated by a decision-focused learning model and directly combine it with a structured covariance estimate. The resulting hybrid matrix is used as an input to mean-variance portfolio optimization. Our focus is on how this particular construction affects portfolio decisions and out-of-sample performance.

\subsection{Data-driven covariance model and decision-focused learning framework}\label{sec:dlcovariance}

Data-driven methods construct covariance matrix from historical returns and other available features. One line of research directly models covariance or semi-covariance matrices \citep{dorta2018structured,monteiro2020stochastic,kim2025estimating}. Another line incorporates data-driven models into structured covariance models by learning components such as factors or loading matrix. For example, NeuralFactors \citep{gopal2024neuralfactors} uses neural networks to learn the loading matrix. In covariance shrinkage estimation, \citet{yang2026neural} use neural networks to determine the shrinkage intensity.

Data-driven covariance models provide flexibility in constructing portfolio risk inputs. By using historical returns and features, they can learn nonlinear mappings from available information to risk estimates and represent ovariance patterns that vary over time. This flexibility, however, also introduces challenges. The estimated covariance matrix depends on the quality and amount of training data, the model specification, and the training objective. Complex models may be sensitive to data noisy or overfitting, particularly when financial data are noisy or market conditions change.

The effectiveness of a data-driven covariance matrix also depends on the portfolio's performance. Conventional ``predict-then-optimize" methods first train a model to predict parameters accurately such as expected returns or covariances and then use the predictions in a separate optimization problem. However, improving prediction accuracy does not necessarily improve the resulting portfolio, because prediction errors can have different consequences for different feasible allocations \citep{elmachtoub2022smart}. Decision-focused learning addresses this distinction by incorporating information from the downstream optimization problem into model training. In portfolio applications, this approach provides a way to learn covariance inputs in relation to the performance of the portfolios they support.

Recent studies have applied decision-focused learning to portfolio problems. In global minimum-variance portfolio construction, \citet{kim2025estimating} use a DLinear-based model to estimate a covariance matrix with a decision-focused training objective. Other studies learn or adjust the spectral representation of a covariance matrix for portfolio construction: \citet{yang2026neural} study neural nonlinear shrinkage, while \citet{bongiorno2026end} consider decision-based cleaning of a correlation matrix. Rather than directly learning a covariance matrix, \citet{anis2025end} use a factor-based risk model in an DFL framework for cardinality-constrained minimum-variance portfolio optimization. These studies illustrate different roles for data-driven models and decision-focused objectives in covariance construction or adjustment. In contrast, our study focuses on the direct matrix-level fusion of two covariance inputs: a structured covariance estimate and a data-driven covariance matrix learned using a portfolio-level objective. The purpose is to examine how their combination affects the risk penalties and portfolio decisions produced by a mean-variance optimization layer.

Several of these studies focus on minimum-variance portfolio construction, where the objective is to reduce the portfolio's risk. Risk control is important, but investors also care about the portfolio's return. Mean-variance optimization incorporates both considerations by combining an expected return term with a variance penalty. It also accommodates investors with different return-risk preferences through adjusting the risk aversion parameter $\lambda$. In this setting, the effect of a covariance matrix depends not only on how it represents asset dependence but also on how its risk penalties interact with the return signal and the investor's risk aversion. Accordingly, we use the hybrid covariance matrix and predicted returns as inputs to a mean–variance portfolio optimization problem and incorporate it as a decision layer within a decision-focused learning framework.

\section{Methodology}\label{sec:methodology}

In this section, we present a decision-focused learning framework that combines structured and data-driven covariance matrices for mean–variance portfolio optimization. The framework consists of a return prediction module, a hybrid covariance prediction module, and a portfolio optimization module. The return prediction module produces expected return estimates, while the hybrid covariance prediction module constructs the structured and data-driven covariance matrices. After scale calibration, the two covariance matrices are combined and used together with the predicted returns in the MV layer to determine portfolio weights. We first present a fixed-weight hybrid scheme and then extend it to an adaptive scheme in which the hybrid weight varies with information available at the decision time.

Let $t=1,\ldots,T$ be the investment periods. At period $t$, the information available before period $t$ is used to construct the input tensor $X_{t-L:t-1} \in \mathbb{R}^{L \times N \times K}$, where $L$ denotes the lookback window length, $N$ is the number of assets, and $K$ is the number of variables associated with each asset. These variables may include historical returns, firm characteristics, and market-state indicators, provided that their values are available at the decision time.

\begin{figure}
    \centering
    \includegraphics[width=0.8\linewidth]{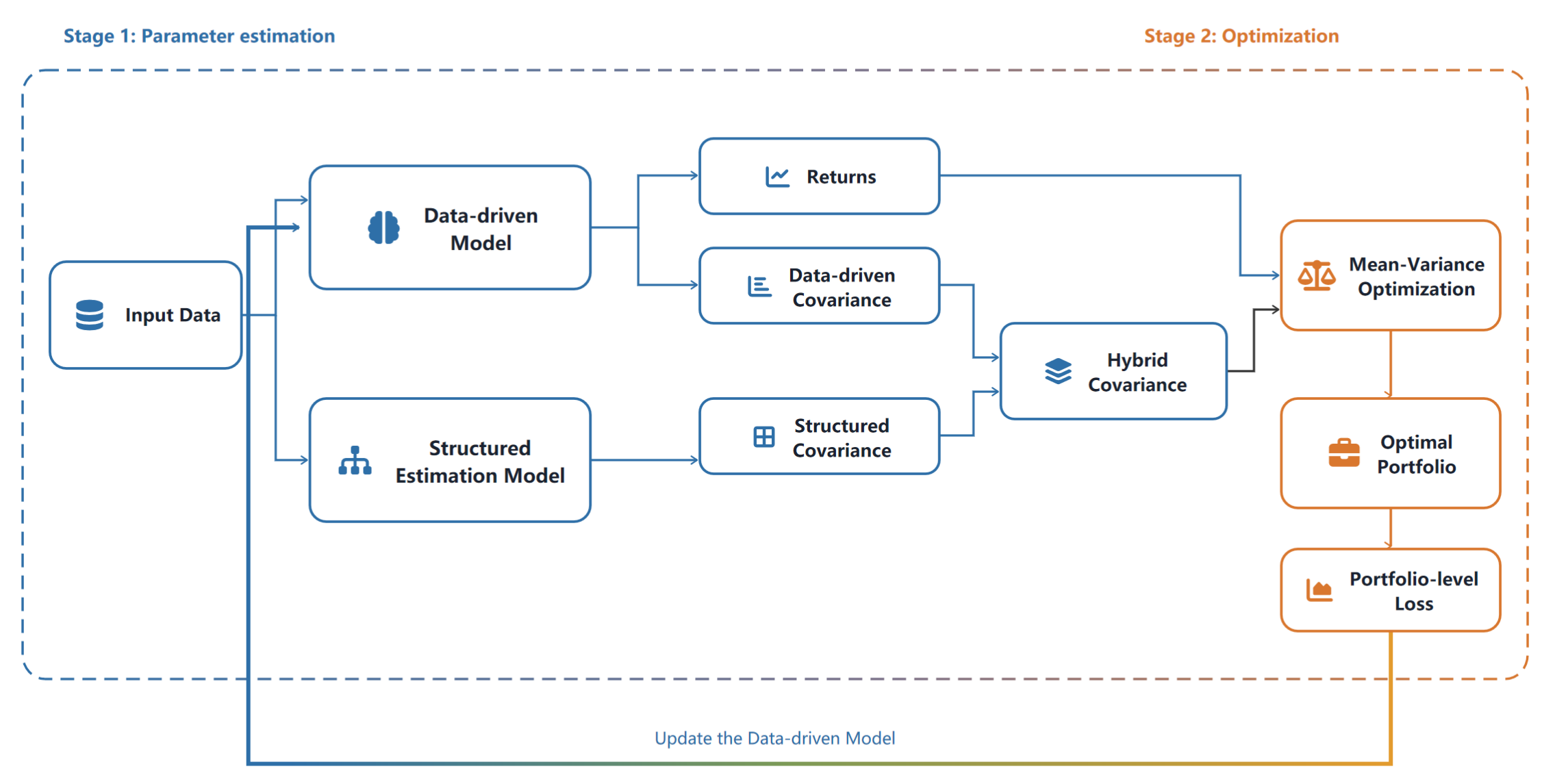}
    \caption{Overview of the decision-focused covariance hybrid framework.}
    \label{fig:framework}
\end{figure}

\subsection{Return prediction module}\label{sec:return_prediction}

The return prediction module characterizes cross-sectional differences in expected returns over the next investment period. Given the input tensor $X_{t-L:t-1}$, the prediction model produces an $N$-dimensional return signal:
\begin{align}
    \widehat{\bfmu}_t = f_{\theta_\mu}\left(X_{t-L:t-1}\right),
\end{align}
where $f_{\theta_\mu}(\cdot)$ the return prediction model with parameters $\theta_\mu$. The predicted return vector and the covariance matrix introduced below refer to the same investment period.

The magnitude of the predicted return signal affects its trade-off with the risk penalty in mean-variance optimization. To control the cross-sectional scale of the signal, we first standardize predicted returns at each decision time:
\begin{align*}
    \widehat{\bfz}_{\mu,t} = \frac{\widehat{\bfmu}_t - \bar{\mu}_t\bfone_N}{s_{\mu,t}},
\end{align*}
where $\bar{\mu}_t$ and $s_{\mu,t}$ are the cross-sectional mean and standard deviation of $\widehat{\bfmu}_t$, respectively. Standardization removes the common cross-sectional level and retains the relative differences among assets.

To express the signal on a scale related to the risk level of the asset universe, we rescale it using the root mean square of the marginal volatilities implied by the structured covariance matrix:
\begin{align*}
    \widetilde{\bfmu}_t = c_{\mu}
    \sqrt{\frac{\operatorname{tr}(\widehat{\Sigma}_t^{\mathrm{str}})}{N}}
    \widehat{\bfz}_{\mu,t},
\end{align*}
where $c_{\mu}>0$ is a scaling parameter and $\widehat{\Sigma}^{\mathrm{str}}_t$ is the structured covariance matrix introduced in Section~\ref{sec:risk_hybrid}. allows the magnitude of the return signal to vary with the overall volatility level of the asset universe. Its influence relative to the risk penalty is jointly determined by $c_{\mu}$, the covariance matrix used in the optimization layer, and the risk aversion parameter.

\subsection{Hybrid covariance prediction module}\label{sec:risk_hybrid}

The hybrid covariance prediction module combines a structured covariance matrix with a data-driven covariance matrix. 
The structured component, denoted as $\widehat{\Sigma}^{\mathrm{str}}$, is obtained from a covariance model that imposes a parsimonious dependence structure. 

The data-driven model constructs a covariance matrix separately from the structured covariance model. Given the input information available at decision time $t$, the model outputs a vector:
\begin{align*}
    \boldsymbol{\ell}_t=f_{\theta_\sigma}(X_{t-L:t-1})\in\mathbb{R}^{N(N+1)/2},
\end{align*}
where $f_{\theta_\sigma}(\cdot)$ is the data-driven model parameterized by $\theta_\sigma$. The vector $\boldsymbol{\ell}_t$ is then arranged into a lower-triangular matrix $L_t\in\mathbb{R}^{N\times N}$.
The data-driven covariance matrix is defined as
\begin{align*}
    \widehat{\Sigma}^{\mathrm{data}}_t=L_tL_t^{\top} + \varepsilon I_N,
\end{align*}
where $I \in \mathbb{R}^{N \times N}$ is the identity matrix and $\varepsilon > 0$ is a numerical regularization parameter. This parameterization ensures that $\widehat{\Sigma}^{\mathrm{data}}_t$ is symmetric positive definite. 

The structured and data-driven covariance matrices are constructed by different input variables and modeling approaches, which can result in different numerical scales. To reduce the influence of scale differences on their hybrid, we calibrate the data-driven matrix to the scale of the structured matrix as follows:
\begin{align*}
    \widetilde{\Sigma}^{\mathrm{data}}_t=\frac{\operatorname{tr}(\widehat{\Sigma}^{\mathrm{str}}_t)}{\operatorname{tr}(\widehat{\Sigma}^{\mathrm{data}}_t)}\widehat{\Sigma}^{\mathrm{data}}_t.
\end{align*}
The calibrated data-driven matrix has the same trace as the structured covariance matrix. This operation removes their difference in total marginal variance while retaining the relative entries of the data-driven matrix.

Therefore, the hybrid covariance matrix is defined as
\begin{align*}
    \widehat{\Sigma}^{\mathrm{hybrid}}_t=\alpha\widetilde{\Sigma}^{\mathrm{data}}_t+(1-\alpha)\widehat{\Sigma}^{\mathrm{str}}_t,\qquad \alpha\in[0,1].
\end{align*}
The parameter $\alpha$ controls the weight assigned to the calibrated data-driven matrix. When $\alpha=0$, the hybrid matrix reduces to a structured covariance matrix; when $\alpha=1$, it reduces to a calibrated data-driven matrix. If the structured covariance matrix is positive semidefinite, the hybrid matrix is also positive semidefinite. 

\subsection{Portfolio optimization module}\label{sec:optimization_training}

Given the calibrated returns and the hybrid covariance matrix, the portfolio weights are obtained from the following long-only mean-variance optimization problem
\begin{align*}
    \bfw^\star_t\in\underset{\bfw \in \W}{\operatorname{argmin}}\left\{-\widetilde{\bfmu}_t^{\top}\bfw + \lambda \bfw^{\top}\widehat{\Sigma}^{\mathrm{hybrid}}_t\bfw\right\},
\end{align*}
where $\lambda>0$ is the risk aversion parameter and the feasible set is
\begin{align*}
    \W=\left\{\bfw\in\R^N:\bfone_N^{\top}\bfw=1,\ \bfw \geq \bf0 \right\}.
\end{align*}
The optimization layer maps the returns and covariance matrix into a fully invested long-only portfolio. A larger value of $\lambda$ places greater relative emphasis on the risk penalty. When $\widehat{\Sigma}^{\mathrm{hybrid}}_t$ is positive semidefinite, the optimization problem is convex. If the hybrid matrix is positive definite, the objective is strongly convex for $\lambda>0$, and the problem has a unique optimal solution over $\W$. 

Let $\bfr_t$ denote the vector of asset returns realized during investment period $t$, and let $r^{\mathrm{free}}_t$ denote the corresponding risk-free return. The realized portfolio excess return is
\begin{align*}
    r^{\mathrm{ex}}_t=\bfr_t^{\top}\bfw^{*}_t-r^{\mathrm{free}}_t.
\end{align*}
The data-driven modules, including the return and data-driven covariance prediction modules, are trained using a portfolio-level objective. Specifically, in this paper, the training loss is the negative Sharpe ratio of realized portfolio excess returns over the training periods:
\begin{align*}
    \L(\theta)=-\frac{\operatorname{mean}_{t\in\mathcal{T}_{\mathrm{train}}}(r^{\mathrm{ex}}_t)}{\operatorname{std}_{t\in\mathcal{T}_{\mathrm{train}}}(r^{\mathrm{ex}}_t)},
\end{align*}
where $\theta=(\theta_\mu,\theta_\sigma)$ collects the parameters of the returns and data-driven covariance prediction models. The loss evaluates the portfolio returns obtained after the optimization layer, rather than directly measuring the prediction error of either model.

Under the fixed-weight scheme, the hybrid parameter $\alpha$ remains fixed, while the return prediction model and the data-driven covariance model are trained jointly. At each period, the two models produce a return signal and a data-driven covariance matrix, respectively. The latter is combined with the structured covariance matrix, and the resulting hybrid matrix and returns are passed to the mean–variance optimization layer to obtain portfolio weights. The parameters of both models are updated using the portfolio-level loss defined above. In this way, the returns and the data-driven covariance matrix are jointly learned according to their combined effect on portfolio performance.

\subsection{Adaptive covariance hybrid extension}\label{subsec:adaptive_alpha}

The fixed-weight scheme uses the same hybrid parameter across decision periods. However, the relative effectiveness of the structured and data-driven covariance matrices may vary with the market conditions and asset dependence patterns evolve. Their effects on portfolio decisions may also differ across levels of risk aversion. We therefore extend the fixed scheme by allowing the hybrid parameter $\alpha$ to vary with features constructed from information available at the decision time.

The adaptive mechanism uses features that summarize differences between the two covariance matrices. First, we measure their structural discrepancy using the normalized Frobenius distance:
\begin{align*}
    d^{\mathrm{cov}}_t=\frac{\left\|\widetilde{\Sigma}^{\mathrm{data}}_t-\widehat{\Sigma}^{\mathrm{str}}_t\right\|_F}{\left\|\widehat{\Sigma}^{\mathrm{str}}_t\right\|_F}.
\end{align*}
A larger value indicates a greater difference between the two covariance matrices relative to the
scale of the structured matrix.

Second, we compare the spectral structures of the covariance matrices using their effective ranks. For a positive semidefinite matrix $\Sigma$, its effective rank is defined as
\begin{align*}
    \operatorname{ER}(\Sigma)=\exp\left(-\sum_{i=1}^{N}p_i\log p_i\right),\qquad p_i=\frac{\lambda_i(\Sigma)}{\sum_{j=1}^{N}\lambda_j(\Sigma)},
\end{align*}
where $\lambda_i(\Sigma)$ is the $i$th eigenvalue and $0\log 0$ is taken to be zero. A higher effective rank indicates a less concentrated distribution of eigenvalues. We define the reference hybrid matrix as
\begin{align*}
    \widetilde{\Sigma}_t^{\mathrm{ref}} = \alpha_{\mathrm{ref}} \widetilde{\Sigma}^{\mathrm{data}}_t + (1-\alpha_{\mathrm{ref}}) \widehat{\Sigma}^{\mathrm{str}}_t,
\end{align*}
where $\alpha_{\mathrm{ref}}$ is fixed. The normalized effective rank difference is
\begin{align*}
    d^{\mathrm{rank}}_t=\frac{\operatorname{ER}(\widetilde{\Sigma}^{\mathrm{ref}}_t)-\operatorname{ER}(\widehat{\Sigma}^{\mathrm{str}}_t)}{N}.
\end{align*}
A positive value indicates that the reference hybrid matrix has a higher effective rank than the structured matrix.

Third, we measure how the portfolio decision changes when the covariance input is switched between the two matrices. We keep the return signal, risk aversion parameter, and constraints unchanged. Let $\bfw^{\mathrm{data}}_t$ and $\bfw^{\mathrm{str}}_t$ be the portfolios obtained from the same mean-variance problem using $\widetilde{\Sigma}^{\mathrm{data}}_t$ and $\widehat{\Sigma}^{\mathrm{str}}_t$, respectively. The weight discrepancy is
\begin{align*}
    d^{w}_t=\left\|\bfw^{\mathrm{data}}_t-\bfw^{\mathrm{str}}_t\right\|_1.
\end{align*}
This feature measures the sensitivity of the portfolio allocation to the choice between the two covariance inputs at the same decision time. 

Denote the state features as:
\begin{align*}
    \bfz_t=\left[d^{\mathrm{cov}}_t,\ d^{\mathrm{rank}}_t,\ d^w_t \right]^{\top}.
\end{align*}
To allow the hybrid weight to account for risk preferences, we further include the risk aversion parameter:
\begin{align*}
    \bfy_t=\left[\bfz_t^{\top},\ \log(1+\lambda)\right]^{\top}.
\end{align*}
All feature variables in $\bfy_t$ are constructed using information available at the decision time.

The adaptive hybrid parameter is generated by a bounded gate:
\begin{align*}
    \alpha_{\phi,t}=\alpha_{\min}+(\alpha_{\max} - \alpha_{\min})\sigma\left(\phi_0+\boldsymbol{\phi}_z^{\top}\bfz_t+\phi_{\lambda}\log(1+\lambda)\right),
\end{align*}
where $\sigma(x)=1/(1+\exp(-x))$, $\phi=(\phi_0,\boldsymbol{\phi}_z,\phi_{\lambda})$ denotes the gate parameters, and $0\leq\alpha_{\min}<\alpha_{\max}\leq 1$.

Using the adaptive weight, the covariance matrix supplied to the mean-variance optimization layer becomes
\begin{align*}
    \widehat{\Sigma}^{\mathrm{hybrid}}_t(\phi)=\alpha_{\phi,t}\widetilde{\Sigma}^{\mathrm{data}}_t+(1-\alpha_{\phi,t})\widehat{\Sigma}^{\mathrm{str}}_t.
\end{align*}
The adaptive scheme retains the same return forecasting module, covariance matrices, scale calibration, and portfolio optimization problem as the fixed-weight scheme. Its difference lies in the determination of the hybrid weight.

We train the adaptive scheme in two stages. In the first stage, the data-driven model for predicting the return and the data-driven covariance matrix are trained jointly under a fixed hybrid weight $\alpha_{\mathrm{base}}$, using the portfolio-level objective described in Section~\ref{sec:optimization_training}. In the second stage, the trained parameters $\theta$ are fixed, and the gate parameters $\phi$ are trained to improve portfolio performance. The corresponding performance criterion over the full training set is the negative Sharpe ratio:
\begin{align*}
    \min_{\phi}\ \L_{\mathrm{adapt}}(\phi;\widehat{\theta}) = -\frac{\operatorname{mean}_{t\in\mathcal{T}_{\mathrm{train}}}(r^{\mathrm{ex}}_t(\phi;\widehat{\theta}))}{\operatorname{std}_{t\in\mathcal{T}_{\mathrm{train}}}(r^{\mathrm{ex}}_t(\phi;\widehat{\theta}))}.
\end{align*}
Here, $\widehat{\theta}$ denotes the model parameters obtained in the first stage. This criterion evaluates the performance of portfolios generated using the adaptive hybrid weight, while the model parameters of the data-driven model remain fixed.

\section{Empirical study}\label{sec:empirical-study}

This section evaluates the proposed decision-focused learning covariance hybrid framework through out-of-sample portfolio experiments and covariance matrix analysis. In the empirical implementation, iTransformer~\citep{liu2024itransformer} is used for both return and data-driven covariance predictions. The structured covariance matrix is constructed using the Fama-French three factors (FF3) model. These choices provide a specific implementation of the modular framework.

The empirical analysis examines four questions. First, how do the proposed fixed $\alpha$ and adaptive hybrid strategies perform relative to traditional and data-driven benchmarks? Second, how does the hybrid covariance matrix affect the portfolio allocations when the predicted returns are fixed? Third, how does the fixed hybrid weight affect portfolio performance across investment styles represented by different levels of risk aversion? Finally, can the adaptive scheme construct competitive portfolios? 

\subsection{Data and experimental setup}\label{data}

The datasets used in this paper are obtained from the Open Source Asset Pricing (OASP) database, the Center for Research in Security Prices (CRSP), and the Fama-French Data Library. The stock daily returns dataset spans from January 4, 2010 to December 31, 2024, covering 3,774 trading days. We select common stocks listed on the New York Stock Exchange using the CRSP conditions sharetype = `NS' and primaryexch = `N'. After merging the stock records with firm characteristics, the resulting dataset contains 3,546 distinct stocks across the full sample period. The number of observed stocks varies over time. 

The original firm-characteristic dataset contains 43 candidate variables covering multiple dimensions, including value, momentum, profitability, volatility, and investor expectations. The characteristics are observed at a monthly frequency and are aligned with daily stock records in the following month. Specifically, characteristics associated with the $m$-th month are matched to daily returns in the $(m+1)$-th month. This alignment introduces a one-month lag between the characteristic month and the return month. Missing characteristic values in the merged daily panel are filled using the cross-sectional median of the corresponding characteristic on each trading day. The 43 characteristics are normalized using daily cross-sectional ranks.

We conduct out-of-sample backtesting using a rolling-window method with 25 windows. Each window consists of 500 trading days for training, followed by 125 trading days for validation and 125 trading days for testing. Within each window, the models are estimated using the training data, model selections are made using the validation data, and performance is evaluated on the test data. For each rolling window, we first identify stocks with no missing return observations in the training period. Among the eligible stocks, we select the Top500 stocks according to the firm-size characteristic to form the investable universe for that window. The selected universe can therefore change across rolling windows. For the selected stocks, if a stock’s return is missing in the validation or test period, its return is set to zero.
We use iTransformer to predict the returns and data-driven covariance matrices. The hidden representation dimension is 16, the number of attention heads is 4, the number of Transformer encoder layers is 1, the dropout rate is 0.3, and the learning rate is $5 \times 10^{-5}$. 
To account for variation across random initializations, the data-driven experiments are repeated using five random seeds. For each seed, portfolio performance metrics are calculated from the corresponding out-of-sample return and weight series. We report the mean of each metric across the five runs.

We consider four investment styles by varying the risk aversion parameter in the mean-variance optimization problem. Specifically, we set $\lambda\in\{10,30,50,100\}$ and refer to the corresponding settings as aggressive, return-dominant, balanced, and risk-dominant, respectively. Under the return and risk scaling used in our framework, a larger $\lambda$ assigns greater relative importance to the risk penalty in the optimization objective. These labels describe the relative return and risk preferences represented by the selected parameter values. For the fixed hybrid weight analysis, we consider $\alpha\in\{0,0.2,0.4,0.6,0.8,1\}$, where $\alpha=0$ uses the structured covariance matrix $\widehat{\Sigma}^{\mathrm{factor}}$ and $\alpha=1$ uses the calibrated data-driven covariance matrix $\widetilde{\Sigma}^{\mathrm{data}}$. Note that, even when $\alpha=1$, the scales of both the return and the data-driven covariance matrix are calibrated using the FF3 covariance matrix as a reference. Thus, this endpoint differs from a purely data-driven approach without factor-based calibration.

We evaluate portfolio performance using annualized return (Ann.Ret.), annualized volatility (Ann.Vol.), Sharpe ratio (Sharpe), maximum drawdown (MaxDD), cumulative return (Cum.Ret.), and average daily turnover (Turnover). We also report performance adjusted for transaction costs. Let $T_{\mathrm{test}}$ denote the number of daily observations in the test period. Starting from initial wealth $W_0=1$, cumulative wealth through test day $t$ is
\begin{align*}
    W_t=\prod_{s=1}^{t}(1+r^p_s), 
\end{align*}
where $r^p_s$ denotes the portfolio return on test day $s$.

The cumulative return and annualized return over the evaluation period are defined as
\begin{align*}
    \operatorname{Cum.Ret.}=W_{T_{\mathrm{test}}}-1,\qquad \operatorname{Ann.Ret.}=W_{T_{\mathrm{test}}}^{252/T_{\mathrm{test}}}-1.
\end{align*}

Annualized volatility and the annualized Sharpe ratio are calculated as
\begin{align*}
    \operatorname{Ann.Vol.}=\sqrt{252}\,s_p,\qquad \operatorname{Sharpe}=\sqrt{252}\frac{\bar r^{\mathrm{ex}}}{s^{\mathrm{ex}}},
\end{align*}
where $s_p$ is the sample standard deviation of daily portfolio returns, and $\bar r^{\mathrm{ex}}$ and $s^{\mathrm{ex}}$ are the sample mean and sample standard deviation of daily excess returns, respectively.

Maximum drawdown measures the largest decline in cumulative wealth from a preceding peak. It is defined as
\begin{align*}
    \operatorname{MaxDD}=\max_{0\leq t\leq T_{\mathrm{test}}}\frac{P_t-W_t}{P_t},\qquad P_t=\max_{0\leq s\leq t}W_s.
\end{align*}

Average daily one-way turnover is defined as
\begin{align*}
    \operatorname{Turnover}=\frac{1}{T_{\mathrm{test}}}\sum_{t=1}^{T_{\mathrm{test}}}\operatorname{TO}_t,\qquad \operatorname{TO}_t=\frac{1}{2}\sum_{i=1}^N\left|w_{i,t}-w^{+}_{i,t-1}\right|,
\end{align*}
where $w_{i,t}$ is the target weight after rebalancing and $w^{+}_{i,t-1}$ is the weight immediately before rebalancing, after accounting for changes in asset values. When the stock universe changes, the summation is taken over the assets relevant to the pre- and post-rebalancing portfolios, with zero weights assigned where an asset is absent.

To assess the effect of trading costs, we report performance under transaction cost rates of 0, 5, and 10 basis points. Let $r^{p,\mathrm{net}}_t(c)$ denote the portfolio return after transaction costs at rate $c$, and define $r^{\mathrm{ex,net}}_t(c)=r^{p,\mathrm{net}}_t(c)-r^{\mathrm{free}}_t$. The transaction-cost-adjusted Sharpe ratio is
\begin{align*}
    \operatorname{NetSharpe}(c)=\sqrt{252}\frac{\bar r^{\mathrm{ex,net}}(c)}{s^{\mathrm{ex,net}}(c)},
\end{align*}
where $\bar r^{\mathrm{ex,net}}(c)$ and $s^{\mathrm{ex,net}}(c)$ are the sample mean and sample standard deviation of daily net excess returns. Net annualized return is calculated from the compounded net return series. 

\subsection{Out-of-sample performance comparison with benchmarks}\label{subsec:oos}

We compare the proposed DF-HCE and Adaptive DF-HCE with traditional and data-driven benchmarks. The traditional benchmarks include the equal-weighted portfolio (EW), the global minimum-variance portfolio (GMV), and mean-variance portfolios constructed using sample (Sample-MV), Ledoit-Wolf shrinkage (Ledoit-Wolf-MV) \citep{ledoit2004honey}, Fama-French three factors (FF3-MV), and Fama-French five factors (FF5-MV) covariance estimates. The data-driven benchmarks include AlphaPortfolio~\citep{cong2021alphaportfolio}, MASTER~\citep{li2024master}, and an iTransformer-based decision-focused learning model (DF-Deep). We also report DF-FF3, which combines iTransformer predicted returns with the FF3 covariance matrix. DF-Deep and DF-FF3 are used for assessing the contribution of hybrid covariance. In this section, we set $\lambda=30$.
All methods use the same dataset and data splits. For each data-driven strategy, experiments are conducted with five random seeds, and performance metrics are reported as averages across the five runs. We also report the sample standard deviation of the Sharpe ratio across seeds. Traditional strategies are deterministic and are evaluated once. MASTER is reproduced using its official implementation, with additional market-state features as required by the model. AlphaPortfolio is independently implemented based on the methodology described in \citet{cong2021alphaportfolio}. 

\begin{figure}[htbp]
\centering
\includegraphics[width=\textwidth]{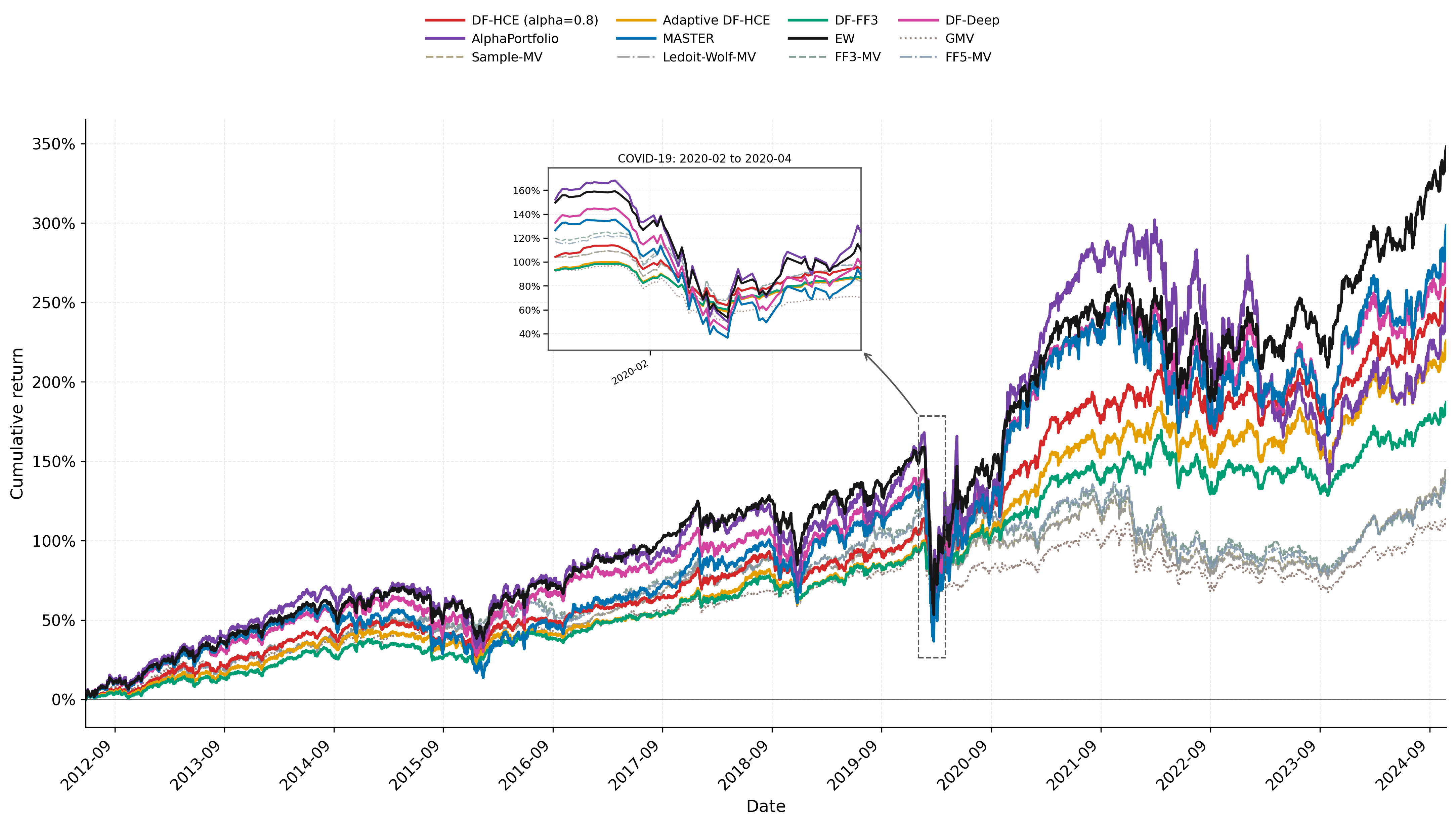}
\caption{Out-of-sample cumulative returns}
\label{fig:benchmark_gross_curves}
\end{figure}

Table~\ref{tab:benchmark_comparison_gross} reports the out-of-sample performance of all strategies over the full test period. 
For DF-HCE, Table~\ref{tab:benchmark_comparison_gross} reports the results with $\alpha=0.8$, which achieves the highest Sharpe ratio among the fixed hybrid weights. The results for other hybrid weights and investment styles are reported in Section~\ref{subsec:lambda_alpha_sensitivity}. At $\lambda=30$, all examined interior hybrid weights achieve higher Sharpe ratios than the traditional and external data-driven benchmarks in Table~\ref{tab:benchmark_comparison_gross}, although not all outperform DF-FF3. 
Figure~\ref{fig:benchmark_gross_curves} shows their cumulative return curves and a closer view of the pandemic period. 
DF-FF3, DF-HCE, and Adaptive DF-HCE achieve higher average Sharpe ratios than the listed benchmarks. Relative to EW and the data-driven benchmarks, they reduce volatility and drawdowns while maintaining competitive returns. Relative to the traditional benchmarks excluding EW, they obtain higher returns by taking on moderately greater risk. The main value of the proposed methods therefore lies in their more competitive return and risk trade-off. Compared with DF-FF3, DF-HCE and Adaptive DF-HCE achieve higher returns by taking on greater risk and turnover.

\begin{table}[htbp]
\centering
\caption{Out-of-sample performance}
\label{tab:benchmark_comparison_gross}
\resizebox{0.9\textwidth}{!}{%
\begin{tabular}{lccccc}
\toprule
Method
& Sharpe (mean \textpm{} std.)
& Ann.Ret. (\%)
& Ann.Vol. (\%)
& MaxDD (\%)
& Turnover (\%) \\
\midrule

\multicolumn{6}{l}{\textit{Traditional benchmark strategies}} \\
EW & 0.7062 & \textcolor{red}{\textbf{12.85}} & 17.40 & 40.77 & \textcolor{red}{\textbf{0.57}} \\
GMV & 0.5543\textsuperscript{*} & 6.34 & \textcolor{blue}{\textbf{9.50}} & 23.60 & \textcolor{blue}{\textbf{1.38}} \\
Sample-MV & 0.6396 & 7.47 & 9.95 & 22.54 & 2.83 \\
Ledoit-Wolf-MV & 0.6323 & 7.45 & 10.07 & 23.08 & 2.89 \\
FF3-MV & 0.5486 & 7.20 & 11.46 & 25.42 & 2.82 \\
FF5-MV & 0.5659 & 7.25 & 11.11 & 24.78 & 2.90 \\

\midrule
\multicolumn{6}{l}{\textit{Data-driven benchmark strategies}} \\
AlphaPortfolio & 0.4990 \textpm{} 0.0907\textsuperscript{**\textdagger\textdagger} & 10.53 & 22.96 & 46.57 & 3.84 \\
MASTER & 0.5741 \textpm{} 0.1715\textsuperscript{**\textdagger} & \textcolor{blue}{\textbf{11.28}} & 19.87 & 42.82 & 10.88 \\
DF-Deep & 0.5937 \textpm{} 0.0438\textsuperscript{**\textdagger} & 11.27 & 18.73 & 41.55 & 9.71 \\

\midrule
\multicolumn{6}{l}{\textit{Proposed strategies}} \\
DF-FF3 & 0.8034 \textpm{} 0.0310 & 8.87 & \textcolor{red}{\textbf{9.48}} & \textcolor{red}{\textbf{19.34}} & 10.69 \\
DF-HCE (alpha = 0.8) & \textcolor{red}{\textbf{0.8491 \textpm{} 0.0612}} & 10.83 & 11.31 & 23.50 & 12.54 \\
Adaptive DF-HCE & \textcolor{blue}{\textbf{0.8249 \textpm{} 0.0405}} & 9.98 & 10.60 & \textcolor{blue}{\textbf{20.83}} & 11.86 \\

\bottomrule
\end{tabular}%
}

\vspace{2mm}
\parbox{0.9\textwidth}{\footnotesize
\textit{Notes:} Superscripts on benchmark Sharpe ratios indicate two-sided HAC-based tests of differences in Sharpe ratios. * and ** denote significance at the 10\% and 5\% levels, respectively, for comparisons with DF-HCE; \textdagger{} and \textdagger\textdagger{} denote the corresponding significance levels for comparisons with Adaptive DF-HCE. 
}
\end{table}

\paragraph{(i) DF-HCE and Adaptive DF-HCE achieve higher out-of-sample Sharpe ratios.}
Table~\ref{tab:benchmark_comparison_gross} shows that the average Sharpe ratios of DF-HCE and Adaptive DF-HCE are 0.8491 and 0.8249, respectively. Both exceed those of the traditional benchmarks (0.5486--0.7062), the data-driven benchmarks (0.4990--0.5937), and DF-FF3 (0.8034). Compared with the traditional benchmarks excluding EW, DF-HCE and Adaptive DF-HCE achieve higher portfolio returns while maintaining a degree of risk control. Compared with the data-driven benchmarks, DF-HCE and Adaptive DF-HCE substantially reduce portfolio volatility and maximum drawdown while preserving competitive returns. Their annualized returns are lower than EW's 12.85\%, but their volatility and maximum drawdown are also markedly lower than EW's 17.40\% and 40.77\%, respectively. Figure~\ref{fig:benchmark_gross_curves} further illustrates the cumulative return paths: EW, MASTER, and DF-Deep achieve relatively high terminal cumulative returns but also experience larger interim declines, whereas the return paths of DF-HCE and Adaptive DF-HCE exhibit comparatively moderate fluctuations.

\paragraph{(ii) Relative to DF-FF3, DF-HCE and Adaptive DF-HCE improve returns while taking on more risk.}
DF-FF3 has an annualized volatility of 9.48\% and a maximum drawdown of 19.34\%, the lowest among the strategies compared. Relative to DF-FF3, DF-HCE and Adaptive DF-HCE increase annualized returns by 1.96 and 1.11 percentage points, respectively. They also raise the Sharpe ratio from 0.8034 to 0.8491 and 0.8249. These improvements are accompanied by higher volatility and maximum drawdowns. Thus, the hybrid covariance models do not uniformly outperform DF-FF3 across all metrics. As shown in Figure~\ref{fig:benchmark_gross_curves}, both DF-HCE and Adaptive DF-HCE end with higher cumulative returns than DF-FF3. This is consistent with the return differences reported in Table~\ref{tab:benchmark_comparison_gross}.

To examine the performance of the proposed methods under a sudden shock and a sustained
market decline, we consider two selected windows: the COVID-19 window from February 1 to
April 30, 2020, and the 2022 market decline from January 1 to October 31, 2022. Table~\ref{tab:market_window_performance} reports the cumulative returns and maximum drawdowns for both windows. In these two selected downturns, DF-HCE and Adaptive DF-HCE exhibit smaller losses and lower maximum drawdowns than the data-driven benchmarks and EW.

\paragraph{(iii) DF-HCE and Adaptive DF-HCE experience smaller losses during the COVID-19 shock and the 2022 market decline.}
The COVID-19 shock and the 2022 market decline represent a sudden disruption and a sustained market decline, respectively, providing complementary settings for evaluating the two methods in downturns.
As shown in Table~\ref{tab:market_window_performance}, during the COVID-19 window, the maximum drawdowns of DF-HCE and Adaptive DF-HCE are 23.50\% and 20.83\%, respectively, below those of the data-driven benchmarks (41.55\%--44.39\%) and EW (40.77\%). Their window losses are 4.73\% and 3.72\%, respectively, also smaller than the losses of these benchmarks (10.19\%--17.27\%). The close-up view of the pandemic period in Figure~\ref{fig:benchmark_gross_curves} illustrates the decline and subsequent rebound of each strategy. DF-FF3 provides even stronger loss control in this window, with a maximum drawdown of 19.34\% and a window loss of 3.34\%.
During the 2022 market decline period, both DF-HCE and Adaptive DF-HCE have lower maximum drawdowns and smaller losses than the data-driven and traditional benchmarks. Their results are broadly comparable to those of DF-FF3: DF-FF3 has a smaller loss, whereas the two DF-HCE strategies have slightly lower maximum drawdowns.

\begin{table}[htbp]
\centering
\caption{Portfolio performance in selected market downturn windows}
\label{tab:market_window_performance}
\setlength{\tabcolsep}{4pt}
\renewcommand{\arraystretch}{0.9}
\resizebox{0.9\textwidth}{!}{%
\begin{tabular}{lcccc}
\toprule
& \multicolumn{2}{c}{COVID-19 window}
& \multicolumn{2}{c}{2022 market decline} \\
\cmidrule(lr){2-3}
\cmidrule(lr){4-5}
Method & Window Ret. (\%) & Max DD (\%) & Window Ret. (\%) & Max DD (\%) \\
\midrule

\multicolumn{5}{l}{\textit{Traditional benchmark strategies}} \\
EW & -15.54 & 40.77 & -12.95 & 20.50 \\
GMV & -11.31 & 23.60 & -15.78 & 20.48 \\
Sample-MV & -9.77 & 21.34 & -16.26 & 21.37 \\
Ledoit-Wolf-MV & -9.86 & 21.53
& -16.57 & 21.73 \\
FF3-MV & -11.49 & 25.42 & -15.60 & 20.89 \\
FF5-MV & -10.47 & 24.11 & -17.25 & 22.08 \\

\midrule
\multicolumn{5}{l}{\textit{Data-driven benchmark strategies}} \\
AlphaPortfolio & -10.19 & 44.39  & -18.65 & 27.00  \\
MASTER & -17.27 & 42.75  & -17.17 & 25.22  \\
DF-Deep & -14.62 & 41.55 & -15.29 & 23.81  \\

\midrule
\multicolumn{5}{l}{\textit{Proposed strategies}} \\
DF-FF3 & \textcolor{red}{\textbf{-3.34}} & \textcolor{red}{\textbf{19.34}} & \textcolor{red}{\textbf{-5.02}} & 14.93 \\
DF-HCE (alpha = 0.8) & -4.73 & 23.50 & \textcolor{blue}{\textbf{-5.35}} & \textcolor{blue}{\textbf{14.42}}  \\
Adaptive DF-HCE & \textcolor{blue}{\textbf{-3.72}} & \textcolor{blue}{\textbf{20.83}} & -5.48 & \textcolor{red}{\textbf{14.36}}  \\

\bottomrule
\end{tabular}%
}
\end{table}

We further assess differences in out-of-sample Sharpe ratios using a two-sided HAC-based test \citep{wright2014a}. The results show that DF-HCE outperforms AlphaPortfolio, MASTER, and DF-Deep at the 5\% significance level, while Adaptive DF-HCE outperforms AlphaPortfolio at the same level. Neither DF-HCE nor Adaptive DF-HCE differs significantly from DF-FF3. 
   
\subsection{Covariance ablation study} \label{subsec:cov-ablation}

To examine how the factor covariance matrix and the iTransformer-based covariance matrix characterize risk and affect portfolio allocation, we hold the network’s calibrated return prediction $\widetilde{\bfmu}$ fixed and vary only the covariance matrix in the optimization objective. We compare portfolios constructed using the sample covariance matrix $\widehat{\Sigma}^{\mathrm{sample}}$, the factor covariance matrix $\widehat{\Sigma}^{\mathrm{factor}}$, the data-driven covariance matrix $\widetilde{\Sigma}^{\mathrm{data}}$, and the hybrid covariance matrix $\widehat{\Sigma}^{\mathrm{hybrid}}$. Denote the portfolios they produce by $\bfw^{\mathrm{sample}}$, $\bfw^{\mathrm{factor}}$, $\bfw^{\mathrm{data}}$, and $\bfw^{\mathrm{hybrid}}$, respectively. 
For notational convenience, let $\Sigma_t^v$ denote the corresponding covariance matrix at time $t$. The cross-risk measures are defined as
\[
q_{v}^k=\frac{1}{T}\sum_{t=1}^T(w_t^k)^\top \Sigma_t^{v} w_t^k, \quad k,v \in \{\mathrm{sample}, \mathrm{factor}, \mathrm{data}, \mathrm{hybrid}\}.
\]
The data-driven risk matrix is scale-calibrated before it is used in portfolio optimization
and cross-risk evaluation. In this section, we set $\alpha=0.8$ and $\lambda=30$. The hybrid covariance matrix is $\widehat{\Sigma}^{\mathrm{hybrid}} = 0.8 \widetilde{\Sigma}^{\mathrm{data}} + 0.2 \widehat{\Sigma}^{\mathrm{factor}}$. The experiments examine three aspects: covariance matrix structure, changes in portfolio weights, and out-of-sample performance. All results in this subsection are based on one random seed (i.e., seed = 1).

\begin{table}[htbp]
    \centering
    \caption{Structural differences between the covariance matrices.}
    \label{tab:risk_structure}
    \begin{threeparttable}
        \small
        \begin{tabular}{lccc}
            \toprule
            Covariance & Diagonal ranking correlation & Off-diagonal covariance correlation & Frobenius distance \\
            \midrule
           $\widehat{\Sigma}^{\mathrm{factor}}$ & 1.0000 & 0.8083 & 0.4098 \\
           $\widetilde{\Sigma}^{\mathrm{data}}$ & 0.0457 & 0.0000 & 1.0028 \\
           $\widehat{\Sigma}^{\mathrm{hybrid}}$ & 0.3366 & 0.4277 & 0.8502 \\
            \bottomrule
        \end{tabular}
        \begin{tablenotes}[flushleft]
            \footnotesize
            \item Benchmark: the sample covariance matrix $\widehat{\Sigma}^{\mathrm{sample}}$.
        \end{tablenotes}
    \end{threeparttable}
\end{table}

\paragraph{(i) The covariance structure represented by $\widetilde{\Sigma}^{\mathrm{data}}$ differs substantially from those represented by $\widehat{\Sigma}^{\mathrm{sample}}$ and $\widehat{\Sigma}^{\mathrm{factor}}$.}

Table~\ref{tab:risk_structure} reports the structural differences between the covariance matrices and $\widehat{\Sigma}^{\mathrm{sample}}$. The results show that $\widehat{\Sigma}^{\mathrm{sample}}$ and $\widehat{\Sigma}^{\mathrm{factor}}$ yield almost identical rankings of individual stock variances, while the correlation between their off-diagonal covariance entries is also high (0.8083). By contrast, the correlation between the individual stock variance rankings of $\widetilde{\Sigma}^{\mathrm{data}}$ and $\widehat{\Sigma}^{\mathrm{sample}}$ is only 0.0457, and the correlation between their off-diagonal covariance entries is close to zero.  This indicates that $\widetilde{\Sigma}^{\mathrm{data}}$ provides a risk representation distinct from those of $\widehat{\Sigma}^{\mathrm{sample}}$ and $\widehat{\Sigma}^{\mathrm{factor}}$.

\begin{table}[htbp]
    \centering
    \caption{Returns and cross-risk measures of portfolios}
    \label{tab:cross_risk_configuration}
    \begin{threeparttable}
        \small
        \begin{tabular}{lrrrr}
            \toprule
            Index & $\bfw^{\mathrm{sample}}$ & $\bfw^{\mathrm{factor}}$ & $\bfw^{\mathrm{data}}$ & $\bfw^{\mathrm{hybrid}}$ \\
            \midrule
            Mean daily return ($10^{-4}$) & 5.0298 & 5.2659 & 7.1451 & 6.1198 \\
            $q_{\mathrm{sample}}^k$ ($10^{-5}$) & 0.3450 & 0.8981 & 16.0684 & 2.8124 \\
            $q_{\mathrm{factor}}^k$ ($10^{-5}$) & 0.5260 & 0.3997 & 15.2469 & 2.1261 \\
            $q_{\mathrm{data}}^k$ ($10^{-5}$) & 2.2930 & 2.1561 & 0.4366 & 0.6286 \\
            $q_{\mathrm{hybrid}}^k$ ($10^{-5}$) & 1.9397 & 1.8049 & 3.3987 & 0.9282 \\
            \bottomrule
        \end{tabular}
    \end{threeparttable}
\end{table}

\paragraph{(ii) $\bfw^{\mathrm{data}}$ has higher returns but receives higher risk estimates under $\widehat{\Sigma}^{\mathrm{sample}}$ and $\widehat{\Sigma}^{\mathrm{factor}}$.}

Table~\ref{tab:cross_risk_configuration} reports the cross-risk evaluation results, revealing how different covariance matrices assess the same portfolio differently. The results show that $q_{\mathrm{data}}^{\mathrm{data}}$ is 4.366$\times 10^{-6}$, substantially lower than $q_{\mathrm{sample}}^{\mathrm{data}}$ (1.6068$\times 10^{-4}$) and $q_{\mathrm{factor}}^{\mathrm{data}}$ (1.5247$\times 10^{-4}$), respectively. In particular, $q_{\mathrm{factor}}^{\mathrm{data}}$ is approximately 34.9 times as large as $q_{\mathrm{data}}^{\mathrm{data}}$. Conversely, $q_{\mathrm{data}}^{\mathrm{factor}}$ is approximately 5.4 times as large as $q_{\mathrm{factor}}^{\mathrm{factor}}$. This bidirectional discrepancy indicates that $\widetilde{\Sigma}^{\mathrm{data}}$ does not uniformly lower risk estimates across all portfolio directions; rather, it changes their relative risk rankings.
This change in risk rankings is accompanied by higher returns. The mean daily return of $\bfw^{\mathrm{data}}$ is 7.1451$\times 10^{-4}$, exceeding the 5.0298$\times 10^{-4}$ and 5.2659$\times 10^{-4}$ achieved by $\bfw^{\mathrm{sample}}$ and $\bfw^{\mathrm{factor}}$, respectively. With predicted returns held fixed, $\bfw^{\mathrm{data}}$ has the highest predicted return among the four portfolios and the lowest $q_{\mathrm{data}}^{\mathrm{data}}$.

\begin{table}[htbp]
    \centering
    \caption{Performance of fixed weights over the future 20 days}
    \label{tab:fixed20_performance}
    \begin{threeparttable}
        \small
        \begin{tabular}{lrrrr}
            \toprule
            Index & $\bfw^{\mathrm{sample}}$ & $\bfw^{\mathrm{factor}}$ & $\bfw^{\mathrm{data}}$ & $\bfw^{\mathrm{hybrid}}$ \\
            \midrule
            Mean excess return ($10^{-4}$) & 1.4641 & 1.4315 & 3.6881 & 3.1659 \\
            Mean variance  & 0.0037 & 0.0038 & 0.0103 & 0.0062 \\
            Sharpe & 0.0878 & 0.0849 & 0.0763 & 0.0945 \\
            \bottomrule
        \end{tabular}
    \end{threeparttable}
\end{table}

\paragraph{(iii) $\bfw^{\mathrm{data}}$ achieves a higher future return, but a lower future Sharpe ratio than $\bfw^{\mathrm{sample}}$ and $\bfw^{\mathrm{factor}}$.}

Table~\ref{tab:fixed20_performance} reports the realized returns and risk over the subsequent 20 days with portfolio weights held fixed. The results show that $\bfw^{\mathrm{data}}$ achieves a higher excess return (3.6881$\times 10^{-4}$) but also has a higher reported variance (0.0103). The increase in returns is insufficient to offset the increase in risk and deliver a higher Sharpe ratio: the Sharpe ratio across the evaluation windows (0.0763) of $\bfw^{\mathrm{data}}$ is lower than those of $\bfw^{\mathrm{sample}}$ (0.0878) and $\bfw^{\mathrm{factor}}$ (0.0849).

\paragraph{(iv) $\bfw^{\mathrm{hybrid}}$ reconstructs the portfolio and improves the return and risk trade-off.}
For any portfolio $\bfw$, we have  
\[
\bfw^\top \widehat{\Sigma}^{\mathrm{hybrid}}\bfw = 0.8 \bfw^\top \widetilde{\Sigma}^{\mathrm{data}}\bfw + 0.2 \bfw^\top \widehat{\Sigma}^{\mathrm{factor}}\bfw.
\]
It follows that
\[
\bfw^\top \widehat{\Sigma}^{\mathrm{hybrid}}\bfw - \bfw^\top \widetilde{\Sigma}^{\mathrm{data}}\bfw = 0.2 \bfw^\top (\widehat{\Sigma}^{\mathrm{factor}} - \widetilde{\Sigma}^{\mathrm{data}})\bfw.
\]

\begin{table}[htbp]
    \centering
    \caption{Out-of-sample performance during the whole test period}
    \label{tab:continuous_performance}
    \begin{threeparttable}
        \small
        \begin{tabular}{lrrrr}
            \toprule
            Index & $\widehat{\Sigma}^{\mathrm{sample}}$ & $\widehat{\Sigma}^{\mathrm{factor}}$ & $\widetilde{\Sigma}^{\mathrm{data}}$ & $\widehat{\Sigma}^{\mathrm{hybrid}}$ \\
            \midrule
            Ann.Ret. (\%) & 4.74 & 5.26 & 8.91 & 9.85 \\
            Ann.Vol. (\%) & 6.43 & 6.53 & 20.20 & 10.39 \\
            MDD (\%) & 13.05 & 14.79 & 44.45 & 18.52 \\
            Sharpe & 0.5462 & 0.6141 & 0.4579 & 0.8280 \\
            NetSharpe (5bps)  & 0.307 & 0.405 & 0.376 & 0.673 \\
            Daily turnover (\%) & 12.18 & 10.84 & 13.18 & 12.82 \\
            \bottomrule
        \end{tabular}
    \end{threeparttable}
\end{table}
Thus, for allocation directions assigned higher risk by $\widehat{\Sigma}^{\mathrm{factor}}$ than $\widetilde{\Sigma}^{\mathrm{data}}$, $\widehat{\Sigma}^{\mathrm{hybrid}}$ increases the risk estimate relative to $\widetilde{\Sigma}^{\mathrm{data}}$ and changes the direction-specific risk penalties in the optimization objective. The results in Table~\ref{tab:cross_risk_configuration} show that, for $\bfw^{\mathrm{data}}$, the risk estimate under $\widehat{\Sigma}^{\mathrm{hybrid}}$ is substantially lower than that under $\widehat{\Sigma}^{\mathrm{factor}}$, but higher than that under $\widetilde{\Sigma}^{\mathrm{data}}$. This indicates that the hybrid approach shifts the allocation away from directions assessed most favorably by $\widetilde{\Sigma}^{\mathrm{data}}$ toward directions more consistent with the risk constraints imposed by $\widehat{\Sigma}^{\mathrm{factor}}$. Moreover, because the risk estimate under $\widehat{\Sigma}^{\mathrm{factor}}$is much higher than that under $\widetilde{\Sigma}^{\mathrm{data}}$, the factor covariance matrix has a substantial influence on the hybrid risk estimate despite its weight of only 20\%.

\begin{figure}
    \centering
    \includegraphics[width=0.8\linewidth]{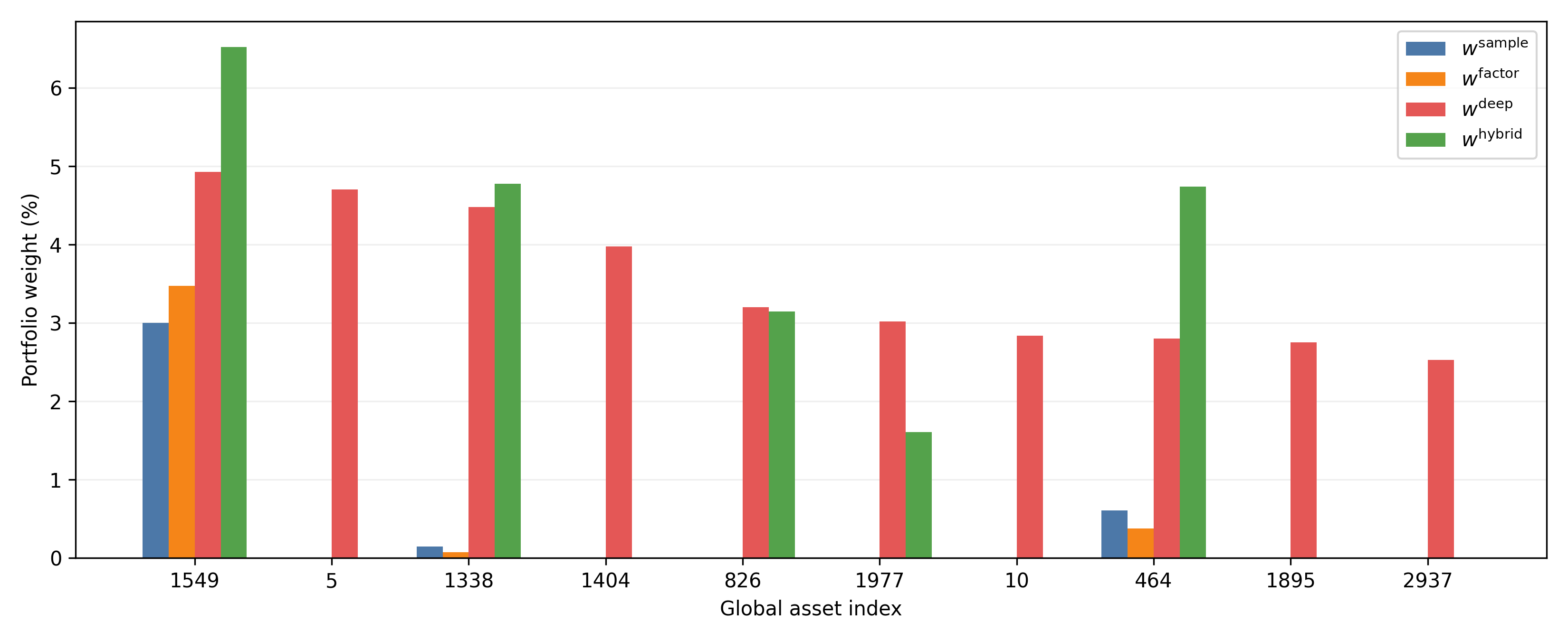}
    \caption{Weights assigned by the four portfolios to the Top10 holdings of $\bfw^{\mathrm{data}}$}
    \label{fig:placeholder}
\end{figure}

Figure~\ref{fig:placeholder} compares the weights assigned by $\bfw^{\mathrm{sample}}$, $\bfw^{\mathrm{factor}}$, and $\bfw^{\mathrm{hybrid}}$ to the Top10 holding of $\bfw^{\mathrm{data}}$. The results show that $\bfw^{\mathrm{hybrid}}$ does not uniformly reduce allocations to the largest holdings of $\bfw^{\mathrm{data}}$; instead, it selectively restructures the portfolio.

Table~\ref{tab:continuous_performance} reports the out-of-sample performance of the different portfolios. Relative to $\bfw^{\mathrm{data}}$, $\bfw^{\mathrm{hybrid}}$ shows improvements in the reported metrics and maintains a higher Sharpe ratio after accounting for transaction costs. Relative to $\bfw^{\mathrm{sample}}$ and $\bfw^{\mathrm{factor}}$, $\bfw^{\mathrm{hybrid}}$ achieves a higher net Sharpe ratio, although it also has higher daily turnover.

\paragraph{(v) The hybrid matrix alleviates numerical ill-conditioning and yields a more dispersed spectrum.} 

\begin{figure}[htbp]
    \centering
    \includegraphics[width=0.9\linewidth]{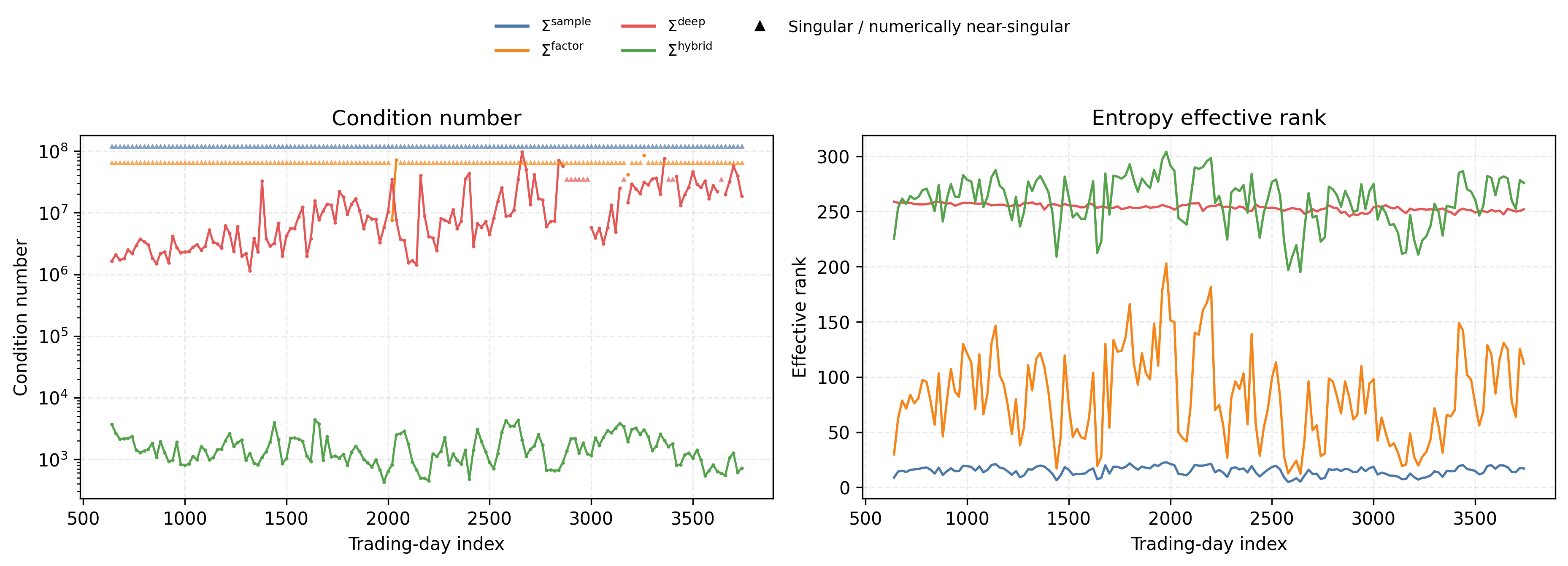}
    \caption{Spectral analysis of the four risk matrices.}
    \label{fig:covariance_spectral}
\end{figure}

Table~\ref{tab:spectral_analysis} shows that $\widehat{\Sigma}^{\mathrm{sample}}$ is structurally singular on all 3,125 evaluation days, whereas $\widehat{\Sigma}^{\mathrm{factor}}$ and $\widetilde{\Sigma}^{\mathrm{data}}$ are classified as numerically near-singular on 3,056 and 212 days, respectively. In contrast, $\widehat{\Sigma}^{\mathrm{hybrid}}$ is neither structurally singular nor classified as near-singular on any evaluation day. Its median condition number is $1.36\times 10^{3}$, substantially lower than those reported for $\widehat{\Sigma}^{\mathrm{factor}}$ and $\widetilde{\Sigma}^{\mathrm{data}}$. $\widehat{\Sigma}^{\mathrm{hybrid}}$ also has the highest mean effective rank (259.97), indicating a more dispersed eigenvalue distribution. These results show that, in the present experiment, the hybrid matrix produces a better-conditioned risk input with a relatively dispersed spectrum; they do not by themselves establish improved portfolio stability or out-of-sample performance. Figure~\ref{fig:covariance_spectral} illustrates how the condition number and effective rank evolve over the test period. The figure samples observations every 20 trading days, so the displayed values and states correspond only to the sampled dates. Together, Table~\ref{tab:spectral_analysis} and Figure~\ref{fig:covariance_spectral} provide aggregate and time-varying perspectives on the spectral properties of the covariance matrices. In this experiment, $\widehat{\Sigma}^{\mathrm{hybrid}}$ provides a better-conditioned covariance input with a relatively dispersed spectrum.

\begin{table}[htbp]
\centering
\caption{Spectral analysis of the four risk matrices}
\label{tab:spectral_analysis}
\scalebox{0.9}{%
\begin{threeparttable}
\begin{tabular}{lrrrr}
\toprule
covariance & Structural singularity & Near singularity & Median Cond.Num & Mean Eff.Rank \\
\midrule
$\widehat{\Sigma}^{\mathrm{sample}}$ & 3125 & 0 & --- & 15.03 \\
$\widehat{\Sigma}^{\mathrm{factor}}$ & 0 & 3056 & \ensuremath{4.15\times 10^{7}} & 84.05 \\
$\widetilde{\Sigma}^{\mathrm{data}}$ & 0 & 212 & \ensuremath{6.54\times 10^{6}} & 253.69 \\
$\widehat{\Sigma}^{\mathrm{hybrid}}$ & 0 & 0 & \ensuremath{1.36\times 10^{3}} & 259.97 \\
\bottomrule
\end{tabular}
\begin{tablenotes}[flushleft]
\footnotesize
\item \textit{Notes:} The results are based on all daily observations. A date is classified as numerically near-singular when $\lambda_{\min}\leq 10^{-8}\lambda_{\max}$. Median condition numbers are calculated only over dates classified as having finite condition numbers. 
\end{tablenotes}
\end{threeparttable}
}
\end{table}

\subsection{Performance of different hybrid weights across investment styles}\label{subsec:lambda_alpha_sensitivity}
\begin{figure}[htbp]
\centering
\includegraphics[width=0.8\textwidth]{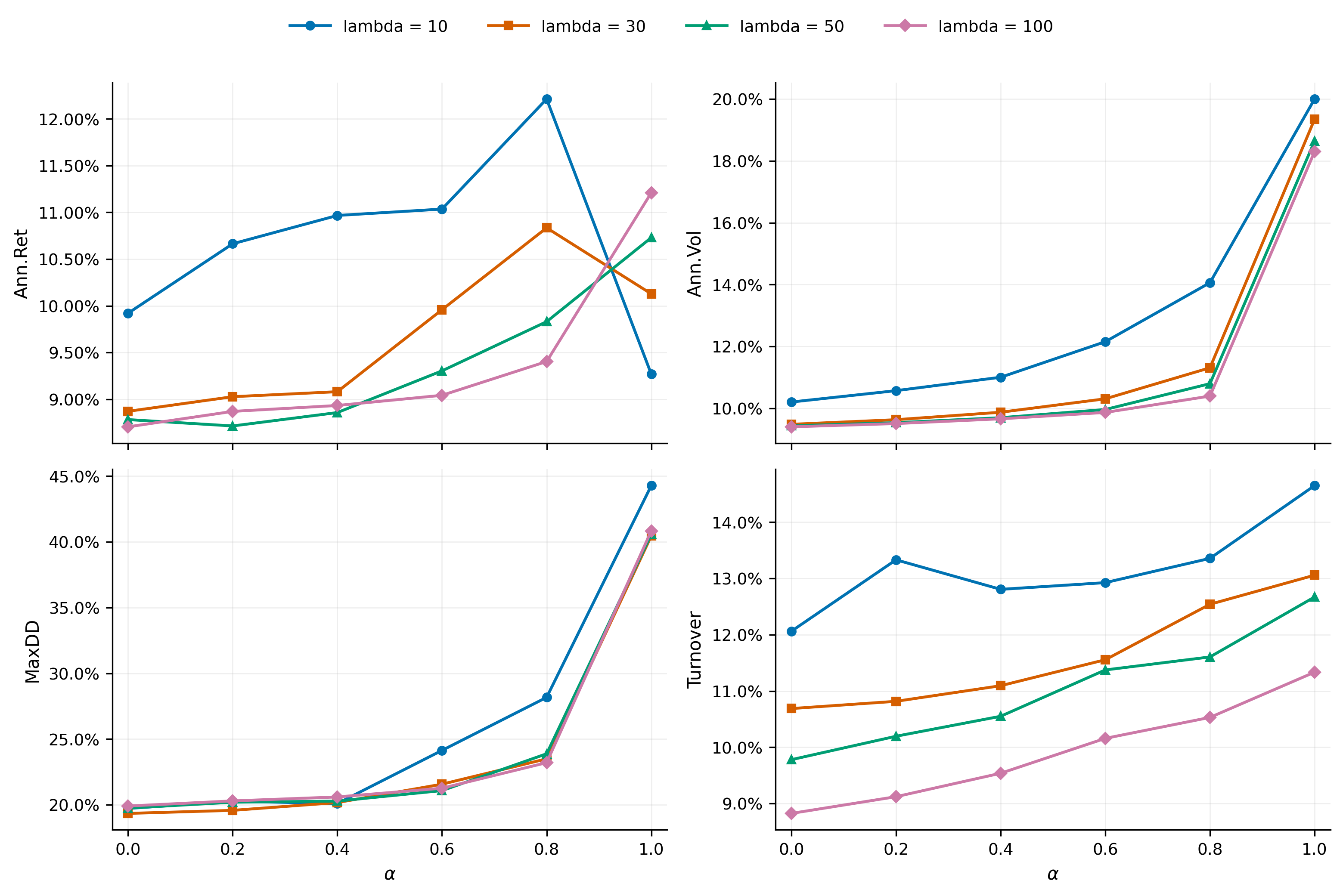}
\caption{Out-of-sample performances across investment styles and hybrid parameters}
\label{fig:lambda_alpha_response}
\end{figure}

This subsection compares the out-of-sample performance of fixed hybrid weights across investment styles to examine how the choice of hybrid weight affects portfolio performance under different risk preferences. We consider $\alpha\in\{0,0.2,0.4,0.6,0.8,1\}$, where $\alpha=0$ and $\alpha=1$ correspond to the FF3 covariance and iTransformer-based covariance, respectively. 
Recall that, when $\alpha=1$, the scales of both the predicted return and the data-driven covariance matrix are calibrated using the FF3 covariance matrix as a reference. Therefore, this endpoint is not a purely data-driven approach without factor-based calibration. For each combination of $\alpha$ and $\lambda$, the model is trained separately. 
Figure~\ref{fig:lambda_alpha_response} shows the annualized returns, annualized volatilities, maximum drawdowns, and average daily one-way turnover.

Overall, interior hybrid weights ($0<\alpha<1$) offer different return-risk trade-offs. Relative to $\alpha=1$, portfolios constructed using the hybrid covariance have lower volatility, maximum drawdown, and turnover; under the aggressive and return-dominant styles, some also achieve higher returns. Relative to $\alpha=0$, at least one interior hybrid weight achieves a higher
annualized return under each investment style, although the magnitude and pattern of improvement depend on the investment style. No single fixed hybrid weight exhibits a consistent overall advantage across all investment styles in the test period. The choice of hybrid weight should therefore reflect the investor's return-risk preferences. This variation motivates the adaptive hybrid weight approach, which allows the weight to respond to decision-related features rather than fixing it in advance.

\paragraph{(i) Relative to $\alpha=1$, hybrid reduces risk and turnover and can improve returns under return-oriented investment styles.}
Across all four investment styles, every interior hybrid weight produces lower annualized volatility, maximum drawdown, and turnover than $\alpha=1$. This indicates that retaining FF3 covariance can reduce the reported risk measures and turnover. Under the aggressive and return-dominant styles ($\lambda\in\{10,30\}$), some interior weights also deliver higher annualized returns. For example, when $\lambda=10$, moving from $\alpha=1$ to $\alpha=0.8$ raises annualized return from 9.27\% to 12.21\%. When $\lambda=30$, the same change raises annualized return from 10.13\% to 10.83\%. Thus, relying exclusively on the data-driven covariance does not yield the highest return under these two styles; retaining some FF3 risk information can improve returns while reducing risk.

By contrast, under the balanced and risk-dominant styles ($\lambda\in\{50,100\}$), $\alpha=1$ yields the highest annualized return, and interior hybrid weights primarily offer a return-risk trade-off. For instance, when $\lambda=100$, moving from $\alpha=1$ to $\alpha=0.8$ reduces annualized return from 11.21\% to 9.40\%, while also substantially reducing volatility and maximum drawdown. Hence, the risk and turnover advantages of interior hybrid relative to the data-driven covariance endpoint are relatively consistent, whereas return improvements depend on the investment style and hybrid weight.

\paragraph{(ii) Relative to $\alpha=0$, incorporating data-driven covariance information generally improves returns but typically increases risk and turnover.}
For each of the four investment styles, at least one interior hybrid weight achieves a higher annualized return than $\alpha=0$. These returns, however, are generally accompanied by increases in volatility, maximum drawdown, and turnover. For example, under the aggressive style ($\lambda=10$), increasing $\alpha$ from 0 to 0.8 raises annualized return from 9.92\% to 12.21\%, while also increasing all three measures. The return improvement from hybrid relative to $\alpha=0$ is therefore not a costless gain: it involves additional risk and trading activity.

\begin{table}[htbp]
\centering
\caption{Gross and net performance across investment styles and hybrid parameters}
\label{tab:alpha_lambda_sensitivity_net}
\setlength{\tabcolsep}{4pt}
\renewcommand{\arraystretch}{0.75}
\begin{adjustbox}{max width=0.8\textwidth}
\begin{tabular}{cccccccc}
\toprule
\multirow{2}{*}{\ensuremath{\lambda}}
& \multirow{2}{*}{\ensuremath{\alpha}}
& \multicolumn{3}{c}{Sharpe}
& \multicolumn{3}{c}{Ann.Ret. (\%)} \\
\cmidrule(lr){3-5}\cmidrule(lr){6-8}
& & 0 bps & 5 bps & 10 bps
& 0 bps & 5 bps & 10 bps \\
\midrule
\multirow{6}{*}{10}
& 0.0 & 0.8478 & 0.6986 & 0.5495
& 9.92 & 8.26 & 6.63 \\
& 0.2 & \textcolor{red}{\textbf{0.8858}} & \textcolor{blue}{\textbf{0.7266}} & \textcolor{blue}{\textbf{0.5674}}
& 10.66 & 8.82 & 7.01 \\
& 0.4 & \textcolor{blue}{\textbf{0.8797}} & \textcolor{red}{\textbf{0.7327}} & \textcolor{red}{\textbf{0.5857}}
& 10.97 & 9.19 & 7.44 \\
& 0.6 & 0.8115 & 0.6770 & 0.5426
& \textcolor{blue}{\textbf{11.03}} & \textcolor{blue}{\textbf{9.24}} & \textcolor{blue}{\textbf{7.48}} \\
& 0.8 & 0.7933 & 0.6732 & 0.5530
& \textcolor{red}{\textbf{12.21}} & \textcolor{red}{\textbf{10.34}} & \textcolor{red}{\textbf{8.50}} \\
& 1.0 & 0.4800 & 0.3871 & 0.2942
& 9.27 & 7.27 & 5.30 \\
\midrule
\multirow{6}{*}{30}
& 0.0 & 0.8034 & 0.6613 & 0.5192
& 8.87 & 7.41 & 5.98 \\
& 0.2 & 0.8075 & 0.6658 & 0.5241
& 9.03 & 7.55 & 6.10 \\
& 0.4 & 0.7950 & 0.6533 & 0.5116
& 9.08 & 7.57 & 6.07 \\
& 0.6 & \textcolor{blue}{\textbf{0.8436}} & \textcolor{blue}{\textbf{0.7021}} & \textcolor{blue}{\textbf{0.5607}}
& 9.96 & \textcolor{blue}{\textbf{8.37}} & \textcolor{blue}{\textbf{6.80}} \\
& 0.8 & \textcolor{red}{\textbf{0.8491}} & \textcolor{red}{\textbf{0.7091}} & \textcolor{red}{\textbf{0.5692}}
& \textcolor{red}{\textbf{10.83}} & \textcolor{red}{\textbf{9.10}} & \textcolor{red}{\textbf{7.39}} \\
& 1.0 & 0.5267 & 0.4417 & 0.3566
& \textcolor{blue}{\textbf{10.13}} & 8.33 & 6.56 \\
\midrule
\multirow{6}{*}{50}
& 0.0 & 0.7987 & \textcolor{red}{\textbf{0.6679}} & \textcolor{red}{\textbf{0.5372}}
& 8.78 & 7.45 & 6.13 \\
& 0.2 & 0.7845 & 0.6496 & 0.5148
& 8.71 & 7.33 & 5.96 \\
& 0.4 & 0.7867 & 0.6495 & 0.5123
& 8.86 & 7.42 & 6.00 \\
& 0.6 & \textcolor{red}{\textbf{0.8093}} & \textcolor{blue}{\textbf{0.6653}} & 0.5213
& 9.30 & 7.75 & 6.21 \\
& 0.8 & \textcolor{blue}{\textbf{0.7994}} & 0.6639 & \textcolor{blue}{\textbf{0.5284}}
& \textcolor{blue}{\textbf{9.83}} & \textcolor{blue}{\textbf{8.24}} & \textcolor{blue}{\textbf{6.67}} \\
& 1.0 & 0.5688 & 0.4830 & 0.3973
& \textcolor{red}{\textbf{10.73}} & \textcolor{red}{\textbf{8.98}} & \textcolor{red}{\textbf{7.25}} \\
\midrule
\multirow{6}{*}{100}
& 0.0 & 0.7935 & \textcolor{blue}{\textbf{0.6751}} & \textcolor{blue}{\textbf{0.5567}}
& 8.70 & 7.50 & 6.31 \\
& 0.2 & \textcolor{red}{\textbf{0.8018}} & \textcolor{red}{\textbf{0.6808}} & \textcolor{red}{\textbf{0.5598}}
& 8.87 & 7.63 & 6.40 \\
& 0.4 & \textcolor{blue}{\textbf{0.7966}} & 0.6720 & 0.5475
& 8.93 & 7.63 & 6.35 \\
& 0.6 & 0.7921 & 0.6623 & 0.5325
& 9.04 & 7.66 & 6.29 \\
& 0.8 & 0.7889 & 0.6611 & 0.5334
& \textcolor{blue}{\textbf{9.40}} & \textcolor{blue}{\textbf{7.96}} & \textcolor{blue}{\textbf{6.54}} \\
& 1.0 & 0.5993 & 0.5213 & 0.4433
& \textcolor{red}{\textbf{11.21}} & \textcolor{red}{\textbf{9.63}} & \textcolor{red}{\textbf{8.08}} \\
\bottomrule
\end{tabular}
\end{adjustbox}
\end{table}

To assess whether the return gains discussed above compensate for the additional risk and trading activity, Table~\ref{tab:alpha_lambda_sensitivity_net} reports Sharpe ratios and annualized returns under different transaction costs. The results show that some hybrid weights provide more favorable Sharpe ratios, and that their advantage persists after transaction costs are taken into account. For example, when $\lambda=30$, $\alpha=0.6$ and $\alpha=0.8$ achieve higher Sharpe ratios than both endpoints at every transaction-cost level examined. Likewise, when $\lambda=10$, $\alpha=0.4$ delivers both a higher Sharpe ratio and a higher annualized return than either endpoint across all examined cost levels. These findings show that selected hybrid weights achieve higher Sharpe ratios than the two endpoints in the examined settings, including after transaction costs.

This advantage is not universal. When $\lambda=50$, $\alpha=0.6$ achieves the highest Sharpe ratio before transaction costs, whereas $\alpha=0$ outperforms all interior hybrid weights at both 5 and 10 bps costs. When $\lambda=100$, $\alpha=0.2$ has a small Sharpe-ratio advantage over $\alpha=0$ at all reported cost levels. Overall, in terms of Sharpe ratio and net Sharpe ratio, hybrid covariance has a clearer advantage under the aggressive and return-dominant styles, whereas DF-FF3 remains competitive under more conservative or cost-sensitive settings.


\subsection{The out-of-sample performance of the adaptive DF-HCE}\label{subsec:adaptive_alpha}
The experiments in Section~\ref{subsec:lambda_alpha_sensitivity} show that out-of-sample performance varies across hybrid weights and that no single fixed weight exhibits a consistent overall advantage across investment styles. This section therefore examines whether adaptive hybrid can construct a competitive portfolio and how the hybrid weight evolves over time. 

\begin{figure}[htbp]
    \centering

    \begin{subfigure}[t]{0.48\textwidth}
        \centering
        \includegraphics[width=\linewidth]
        {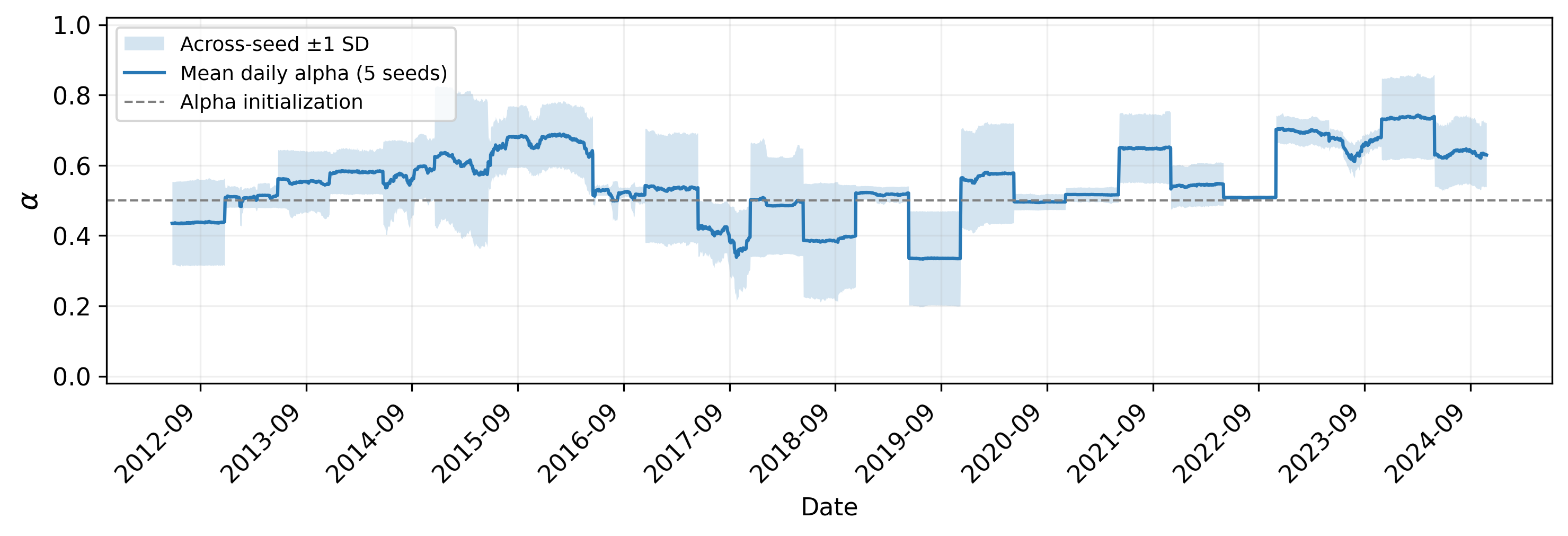}
        \caption{$\lambda = 10$}
        \label{fig:alpha_lambda1}
    \end{subfigure}
    \hfill
    \begin{subfigure}[t]{0.48\textwidth}
        \centering
        \includegraphics[width=\linewidth]
        {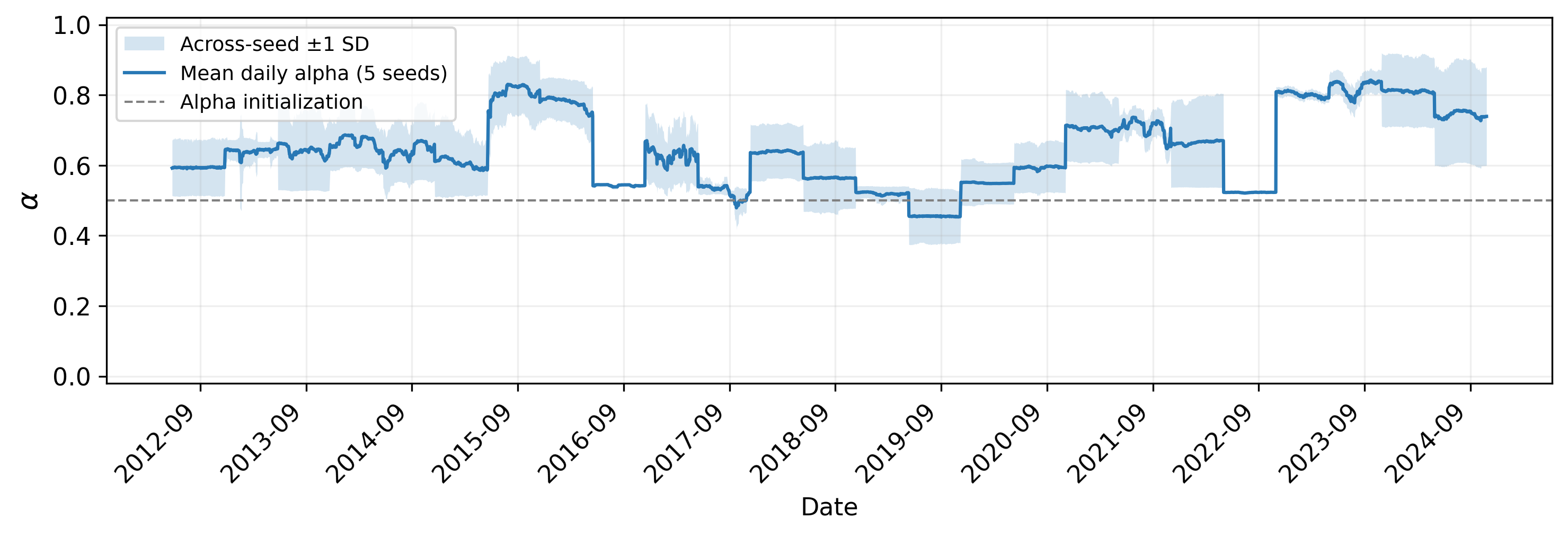}
        \caption{$\lambda = 30$}
        \label{fig:alpha_lambda10}
    \end{subfigure}

    \par\vspace{0.1cm}

    \begin{subfigure}[t]{0.48\textwidth}
        \centering
        \includegraphics[width=\linewidth]
        {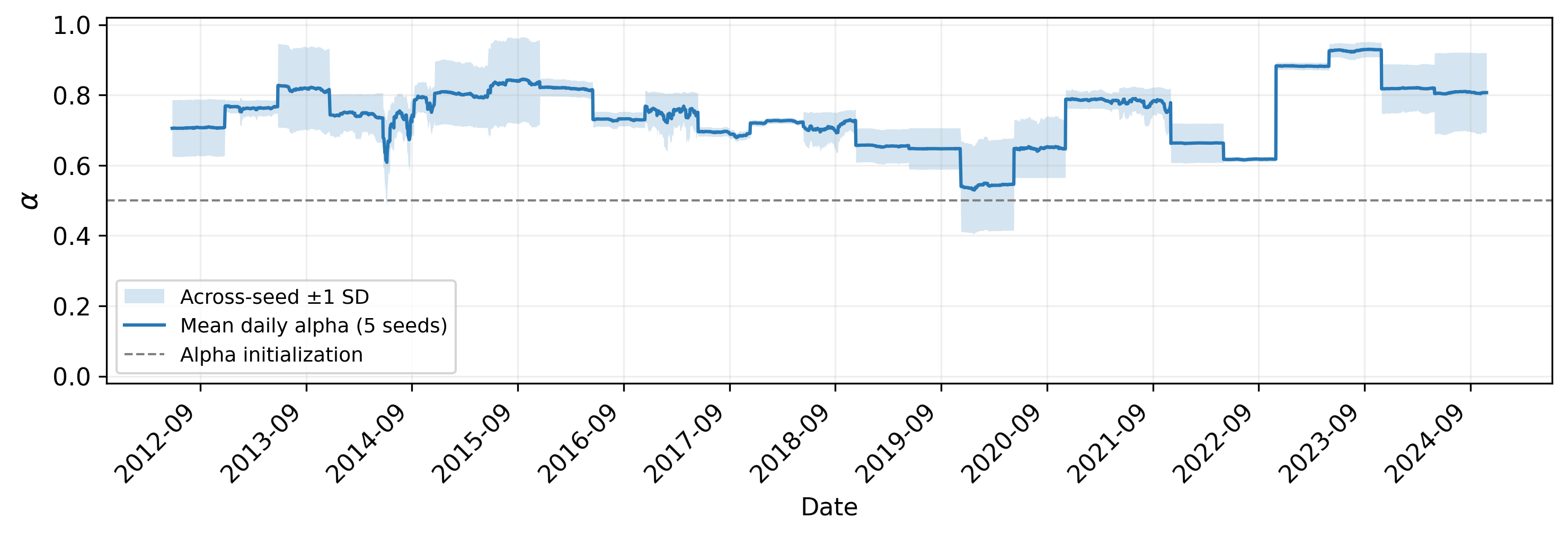}
        \caption{$\lambda = 50$}
        \label{fig:alpha_lambda50}
    \end{subfigure}
    \hfill
    \begin{subfigure}[t]{0.48\textwidth}
        \centering
        \includegraphics[width=\linewidth]
        {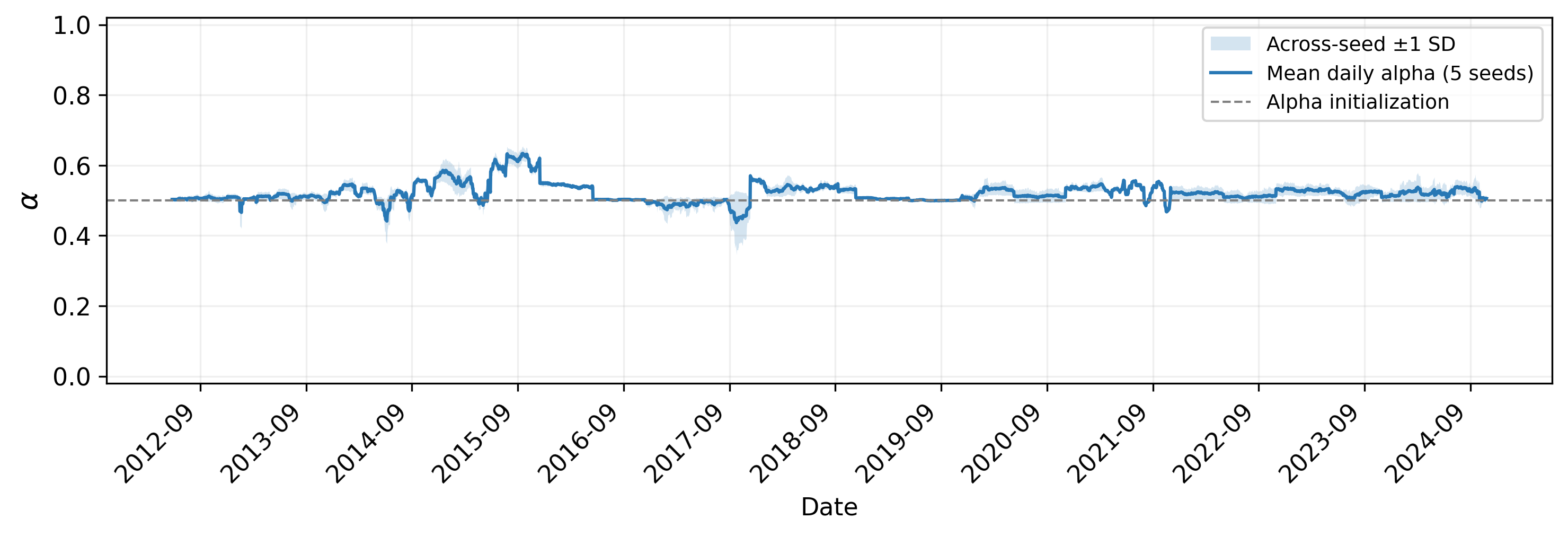}
        \caption{$\lambda = 100$}
        \label{fig:alpha_lambda100}
    \end{subfigure}

    \caption{Time series of the adaptive mixing coefficient
    under different regularization parameters.}
    \label{fig:adaptive_alpha_comparison}
\end{figure}

Table~\ref{tab:adaptive_alpha_performance} reports the out-of-sample performance of adaptive DF-HCE, with net performance calculated using a transaction cost of 5 bps. Figure~\ref{fig:adaptive_alpha_comparison} shows the evolution of the mean hybrid weight across seeds under different investment styles. Adaptive DF-HCE achieves higher annualized returns than DF-FF3 under all four styles, and this return advantage persists after accounting for the specified transaction cost. Its Sharpe ratio also exceeds that of DF-FF3 when $\lambda\in\{30,50,100\}$, and exceeds that of the data-driven risk-matrix endpoint ($\alpha=1$) under all four styles. 
The plotted mean weights remain between zero and one, although their levels and time variation differ across the four settings. Thus, adaptive DF-HCE offers a competitive alternative to portfolios based on either risk-matrix endpoint.

\paragraph{(i) Adaptive DF-HCE increases annualized returns relative to DF-FF3 under all four investment styles, with Sharpe-ratio gains under most styles.}
Relative to DF-FF3, adaptive DF-HCE increases annualized returns by 0.65, 1.11, 1.02, and 0.33 percentage points across the four investment styles, respectively. These return advantages persist after accounting for a 5 bps transaction cost. When $\lambda\in\{30,50,100\}$, the Sharpe ratios rise from 0.8034, 0.7987, and 0.7935 for DF-FF3 to 0.8249, 0.8029, and 0.7956 for adaptive DF-HCE, respectively. 
After transaction costs, however, adaptive DF-HCE retains a net Sharpe ratio advantage only when $\lambda=30$. 
Compared with DF-HCE, the Sharpe ratio of adaptive DF-HCE is close to the best value when $\lambda\in\{50,100\}$. When $\lambda=10$, its Sharpe ratio falls below that of DF-FF3 and some DF-HCE strategies. Therefore, adaptive DF-HCE can produce competitive portfolios, but higher returns do not necessarily translate into better risk-adjusted performance after costs.

\paragraph{(ii) The adaptive weights remain between the two endpoints and tend to favor the data-driven covariance.}
The mean hybrid weights under the four investment styles are 0.5512, 0.6468, 0.7443, and 0.5213, respectively. All exceed 0.5, with the strongest tilt toward the data-driven covariance occurring when $\lambda\in\{30,50\}$. 
Figure~\ref{fig:adaptive_alpha_comparison} shows that the mean weights vary over time. Throughout the observed period, the adaptive method retains information from both covariance matrices while varying their relative contributions.

\begin{table}[htbp]
\centering
\caption{Out-of-sample performance and mean hybrid parameters of adaptive DF-HCE}
\label{tab:adaptive_alpha_performance}
\small
\resizebox{0.9\textwidth}{!}{%
\begin{tabular}{cccccccccc}
\toprule
$\lambda$
& Sharpe
& Ann.Ret(\%)
& Ann.Vol(\%)
& Turnover
& MaxDD
& Mean $\alpha$
& NetSharpe
& NetRet(\%) \\
\midrule
10  & 0.8058 & 10.57 & 11.65 & 0.1274 & 0.2056 & 0.5512 & 0.6678 & 8.81 \\
30  & 0.8249 &  9.98 & 10.60 & 0.1186 & 0.2083 & 0.6468 & 0.6839 & 8.35 \\
50  & 0.8029 &  9.80 & 10.71 & 0.1220 & 0.2129 & 0.7443 & 0.6593 & 8.13 \\
100 & 0.7956 &  9.03 &  9.80 & 0.1013 & 0.2112 & 0.5213 & 0.6652 & 7.65 \\
\bottomrule
\end{tabular}%
}
\end{table}

\section{Conclusion}\label{sec:conclusion}

In this paper, we propose a decision-focused covariance hybrid framework for mean-variance portfolio optimization. The framework directly combines a structured covariance estimate with a data-driven covariance matrix learned using a portfolio-level objective. It brings together structured covariance information and a decision-focused data-driven covariance representation in a common covariance input for portfolio construction. We further develop an adaptive extension in which we adjust the hybrid weight by the decision features. 
Empirical results show that, across investment styles, DF-HCE achieves higher out-of-sample Sharpe ratios than the benchmarks. The Adaptive DF-HCE also delivers competitive risk-adjusted performance, offering a flexible alternative to the fixed hybrid weights. Further analysis shows that the two component matrices provide distinct representations of risk. Their combination changes the portfolio allocations and risk penalties. Spectral analysis further indicates that the hybrid matrix exhibits better numerical conditioning.
Future research could incorporate transaction costs into training and develop a theoretical basis for selecting and refining the features used in the adaptive scheme.

\bibliographystyle{informs2014} 
\bibliography{references}

@article{broadie1993computing,
  title={Computing efficient frontiers using estimated parameters},
  author={Broadie, Mark},
  journal={Annals of operations research},
  volume={45},
  number={1},
  pages={21-58},
  year={1993},
  publisher={Springer}
}

@article{wright2014a,
  title={A test for the equality of multiple Sharpe ratios},
  author={Wright, John Alexander and Yam, Sheung Chi Phillip and Yung, Siu Pang},
  journal={The Journal of Risk},
  year={2014},
  publisher={Incisive Media. The Journal's web site is located at http://www~…}
}

@inproceedings{liu2024itransformer,
  title={{iTransformer}: Inverted transformers are effective for time series forecasting},
  author={Liu, Yong and Hu, Tengge and Zhang, Haoran and Wu, Haixu and Wang, Shiyu and Ma, Lintao and Long, Mingsheng},
  booktitle={International conference on learning representations},
  volume={2024},
  pages={11116--11140},
  year={2024}
}

@inproceedings{gopal2024neuralfactors,
  title={Neuralfactors: A novel factor learning approach to generative modeling of equities},
  author={Gopal, Achintya},
  booktitle={Proceedings of the 5th ACM International Conference on AI in Finance},
  pages={99--107},
  year={2024}
}

@article{lopez2016building,
  title={Building diversified portfolios that outperform out-of-sample},
  author={Lopez de Prado, Marcos},
  journal={Journal of Portfolio Management},
  year={2016}
}

@article{zou2006sparse,
  title={Sparse principal component analysis},
  author={Zou, Hui and Hastie, Trevor and Tibshirani, Robert},
  journal={Journal of computational and graphical statistics},
  volume={15},
  number={2},
  pages={265--286},
  year={2006},
  publisher={Taylor \& Francis}
}

@article{laloux2000random,
  title={Random matrix theory and financial correlations},
  author={Laloux, Laurent and Cizeau, Pierre and Potters, Marc and Bouchaud, Jean-Philippe},
  journal={International Journal of Theoretical and Applied Finance},
  volume={3},
  number={03},
  pages={391--397},
  year={2000},
  publisher={World Scientific}
}

@article{kan2007optimal,
  title={Optimal portfolio choice with parameter uncertainty},
  author={Kan, Raymond and Zhou, Guofu},
  journal={Journal of Financial and Quantitative Analysis},
  volume={42},
  number={3},
  pages={621-656},
  year={2007},
  publisher={Cambridge University Press}
}

@article{fama1993common,
  title={Common risk factors in the returns on stocks and bonds},
  author={Fama, Eugene F and French, Kenneth R},
  journal={Journal of financial economics},
  volume={33},
  number={1},
  pages={3-56},
  year={1993},
  publisher={Elsevier}
}

@article{fan2008high,
  title={High dimensional covariance matrix estimation using a factor model},
  author={Fan, Jianqing and Fan, Yingying and Lv, Jinchi},
  journal={Journal of Econometrics},
  volume={147},
  number={1},
  pages={186-197},
  year={2008},
  publisher={Elsevier}
}

@article{de2021factor,
  title={Factor models for portfolio selection in large dimensions: The good, the better and the ugly},
  author={De Nard, Gianluca and Ledoit, Olivier and Wolf, Michael},
  journal={Journal of Financial Econometrics},
  volume={19},
  number={2},
  pages={236--257},
  year={2021},
  publisher={Oxford University Press}
}

@article{gu2021autoencoder,
  title={Autoencoder asset pricing models},
  author={Gu, Shihao and Kelly, Bryan and Xiu, Dacheng},
  journal={Journal of Econometrics},
  volume={222},
  number={1},
  pages={429--450},
  year={2021},
  publisher={Elsevier}
}

@article{monteiro2020stochastic,
  title={Stochastic segmentation networks: Modelling spatially correlated aleatoric uncertainty},
  author={Monteiro, Miguel and Le Folgoc, Lo{\"\i}c and Coelho de Castro, Daniel and Pawlowski, Nick and Marques, Bernardo and Kamnitsas, Konstantinos and Van der Wilk, Mark and Glocker, Ben},
  journal={Advances in neural information processing systems},
  volume={33},
  pages={12756--12767},
  year={2020}
}

@inproceedings{dorta2018structured,
  title={Structured uncertainty prediction networks},
  author={Dorta, Garoe and Vicente, Sara and Agapito, Lourdes and Campbell, Neill DF and Simpson, Ivor},
  booktitle={2018 IEEE/CVF Conference on Computer Vision and Pattern Recognition},
  pages={5477--5485},
  year={2018},
  organization={IEEE}
}

@article{kim2018forecasting,
  title={Forecasting the volatility of stock price index: A hybrid model integrating LSTM with multiple GARCH-type models},
  author={Kim, Ha Young and Won, Chang Hyun},
  journal={Expert Systems with Applications},
  volume={103},
  pages={25--37},
  year={2018},
  publisher={Elsevier}
}

@inproceedings{li2024master,
  title={{MASTER}: Market-guided stock transformer for stock price forecasting},
  author={Li, Tong and Liu, Zhaoyang and Shen, Yanyan and Wang, Xue and Chen, Haokun and Huang, Sen},
  booktitle={Proceedings of the AAAI conference on artificial intelligence},
  volume={38},
  number={1},
  pages={162--170},
  year={2024}
}

@article{cong2021alphaportfolio,
  title={{AlphaPortfolio}: Direct construction through deep reinforcement learning and interpretable {AI}},
  author={Cong, Lin William and Tang, Ke and Wang, Jingyuan and Zhang, Yang},
  journal={Available at SSRN 3554486},
  year={2021}
}

@inproceedings{kim2025estimating,
  title={Estimating covariance for global minimum variance portfolio: A decision-focused learning approach},
  author={Kim, Juchan and Tae, Inwoo and Lee, Yongjae},
  booktitle={Proceedings of the 6th ACM International Conference on AI in Finance},
  pages={105--113},
  year={2025}
}

@article{conlon2025asset,
  title={Asset allocation with factor-based covariance matrices},
  author={Conlon, Thomas and Cotter, John and Kynigakis, Iason},
  journal={European Journal of Operational Research},
  volume={325},
  number={1},
  pages={189--203},
  year={2025},
  publisher={Elsevier}
}

@article{ledoit2012nonlinear,
  title={Nonlinear shrinkage estimation of large-dimensional covariance matrices},
  author={Ledoit, Olivier and Wolf, Michael},
  journal={The Annals of Statistics },
  volume={40},
  number={2},
  pages={1024--1060},
  year={2012}
}

@article{kelly2019characteristics,
  title={Characteristics are covariances: A unified model of risk and return},
  author={Kelly, Bryan T and Pruitt, Seth and Su, Yinan},
  journal={Journal of Financial Economics},
  volume={134},
  number={3},
  pages={501--524},
  year={2019},
  publisher={Elsevier}
}

@article{zhu2024dynamic,
  title={Dynamic {ETF} Portfolio Optimization Using enhanced {Transformer}-Based Models for Covariance and Semi-Covariance Prediction},
  author={Zhu, Jiahao and Wu, Hengzhi},
  journal={arXiv preprint arXiv:2411.19649},
  year={2024}
}

@article{yang2026neural,
  title={Neural Nonlinear Shrinkage of Covariance Matrices for Minimum Variance Portfolio Optimization},
  author={Yang, Liusha and Zhao, Siqi and Chai, Shuqi},
  journal={arXiv preprint arXiv:2601.15597},
  year={2026}
}

@article{elmachtoub2022smart,
  title={{Smart “Predict, then Optimize”}},
  author={Elmachtoub, Adam N and Grigas, Paul},
  journal={Management Science},
  volume={68},
  number={1},
  pages={9--26},
  year={2022},
  publisher={INFORMS}
}

@article{fama2015five,
  title={A five-factor asset pricing model},
  author={Fama, Eugene F and French, Kenneth R},
  journal={Journal of financial economics},
  volume={116},
  number={1},
  pages={1--22},
  year={2015},
  publisher={Elsevier}
}

@article{carhart1997persistence,
  title={On persistence in mutual fund performance},
  author={Carhart, Mark M},
  journal={The Journal of finance},
  volume={52},
  number={1},
  pages={57--82},
  year={1997},
  publisher={Wiley Online Library}
}

@article{anis2025end,
  title={End-to-end, decision-based, cardinality-constrained portfolio optimization},
  author={Anis, Hassan T and Kwon, Roy H},
  journal={European Journal of Operational Research},
  volume={320},
  number={3},
  pages={739--753},
  year={2025},
  publisher={Elsevier}
}

@article{ardia2017the,
  title={The impact of covariance misspecification in risk-based portfolios},
  author={Ardia, David and Bolliger, Guido and Boudt, Kris and Gagnon-Fleury, Jean-Philippe},
  journal={Annals of Operations Research},
  volume={254},
  number={1},
  pages={1--16},
  year={2017},
  publisher={Springer}
}

@article{markowitz1952,
	title={Portfolio selection},
	author={Harry Markowitz},
	journal={The Journal of Finance},
	volume={7},
	pages={77-91},
	year={1952},
	publisher={Wiley Online Library}
}

@article{ledoit2004honey,
	title="Honey, I shrunk the sample covariance matrix",
	author="Olivier Ledoit and Michael Wolf",
	journal="The Journal of Portfolio Management",
	volume="30",
	number="4",
	pages="110--119",
	notes="Sourced from Microsoft Academic - https://academic.microsoft.com/paper/3124158341",
	year="2004"
}

@article{demiguel2009a,
	title="A Generalized Approach to Portfolio Optimization: Improving Performance by Constraining Portfolio Norms",
	author="Victor {DeMiguel} and Lorenzo {Garlappi} and Francisco J. {Nogales} and Raman {Uppal}",
	journal="Management Science",
	volume="55",
	number="5",
	pages="798--812",
	notes="Sourced from Microsoft Academic - https://academic.microsoft.com/paper/2127000676",
	year="2009"
}

@article{Fan2012Vast,
  title = {Vast Portfolio Selection With Gross-Exposure Constraints},
  author = {Jianqing Fan and Jingjin Zhang and Ke Yu}, 
  journal = {Journal of the American Statistical Association}, 
  number = {498}, 
  pages = {592--606},
  publisher = {[American Statistical Association, Taylor & Francis, Ltd.]},
  volume = {107},
  year = {2012}
}

@article{bongiorno2026end,
  title={End-to-end large portfolio optimization for variance minimization with neural networks through covariance cleaning},
  author={Bongiorno, Christian and Manolakis, Efstratios and Mantegna, Rosario Nunzio},
  journal={The Journal of Finance and Data Science},
  pages={100179},
  year={2026},
  publisher={Elsevier}
}

\end{document}